\documentclass[12pt, reqno]{amsart}
\usepackage{preemble}
\usepackage{microtype}
\begin{document}

\title[Noise sensitivity on affine Weyl groups]
      {Noise sensitivity on affine Weyl groups}
\author{Ryokichi Tanaka}
\address{Department of Mathematics, Kyoto University, Kyoto 606-8502 JAPAN}
\email{rtanaka@math.kyoto-u.ac.jp}
\date{\today}

\maketitle

\begin{abstract}
We show that on every affine Weyl group natural random walks are noise sensitive in total variation.
\end{abstract}

\section{Introduction}

Let $\Gamma$ be a countable group and $\mu$ be a probability measure on it.
A $\mu$-{\bf random walk} $\{w_n\}_{n \in \Z_+}$ starting from the identity $\id$ is defined by $w_n:=x_1 \cdots x_n$ and $w_0:=\id$ for an independent, identically distributed sequence $x_1, x_2, \dots$ with the common law $\mu$.
The distribution of $w_n$ is the $n$-fold convolution $\mu_n:=\mu^{\ast n}$.
The noise sensitivity problem for random walks on groups asks:
Does resampling a small fraction of increments $x_1, x_2, \dots$ produce an almost independent copy of $w_n$ or a highly correlated copy of $w_n$?

The precise definition is as follows.
For a real $\rho \in [0, 1]$,
let 
\[
\pi^\rho:=\rho (\mu \times \mu) +(1-\rho)\mu_{\diag} \quad \text{on $\Gamma \times \Gamma$},
\]
where $\mu \times \mu$ denotes the product measure and $\mu_{\diag}((x, y)):=\mu(x)$ if $x=y$ and $0$ otherwise.
Let us consider a $\pi^\rho$-random walk $\{\wb_n\}_{n\in \Z_+}$ starting from the identity on $\Gamma \times \Gamma$.
We say that the $\mu$-random walk is {\bf noise sensitive} in total variation if
\[
\lim_{n \to \infty}\|\pi^\rho_n-\mu_n \times \mu_n\|_\TV =0 \quad \text{for all $\rho \in (0, 1]$}.
\]
In the above, the total variation distance coincides with the half of $\ell^1$-norm,
\[
\|\nu_1-\nu_2\|_\TV=\max_{A \subset \Gamma\times \Gamma}|\nu_1(A)-\nu_2(A)|=\frac{1}{2}\sum_{\xb \in \Gamma\times \Gamma}|\nu_1(\xb)-\nu_2(\xb)|, 
\]
where $\nu_i$ are probability measures on $\Gamma\times \Gamma$, $i=1, 2$.
Since we use only total variation distance in the definitions, let us simply say noise sensitive if there is no danger of confusion.

If a $\mu$-random walk on $\Gamma$ is noise sensitive,
then informally speaking, the situation is as in the following: For each fixed $\rho \in (0, 1]$ even though it is close to $0$,
resampling a $\rho$-portion of increments produces an asymptotically independent copy of the original $\mu$-random walk.

The definition of noise sensitivity for random walks on groups was introduced by Benjamini and Brieussel in \cite[Definition 2.1]{BenjaminiBrieussel}.
In their paper, they discuss the notion and the variants not only in total variation (there it is called $\ell^1$-noise sensitivity) but also in other distances or in terms of entropy.
See also the related discussion in \cite[Section 3.3.4]{KalaiICM2018}.
It has been observed that on finite groups random walks are noise sensitive in all natural definitions \cite[Proposition 5.1]{BenjaminiBrieussel}.
This raises a challenge to find finitely generated {\it infinite} groups on which random walks are noise sensitive.
Further the problem becomes more restrictive by measuring the distance in total variation.
A simple observation using the central limit theorem shows that standard random walks on finite rank free abelian groups $\Z^m$ are {\it not} noise sensitive \cite[Theorem 1.1 (1)]{BenjaminiBrieussel}.
Benjamini and Brieussel have shown that on the infinite dihedral group some lazy simple random walk is noise sensitive \cite[Theorem 1.4]{BenjaminiBrieussel}.
So far, this has been the only known random walk which is noise sensitive in total variation on a finitely generated infinite group.
We provide a class of such groups on which natural random walks are noise sensitive.

\begin{theorem}\label{Thm:NSonWeylmain}
Let $\Gamma$ be an affine Weyl group with the standard set of generators $S$, and $\mu$ be a probability measure on $\Gamma$ such that
the support of $\mu$ equals $S\cup\{\id\}$.
For all $\rho \in (0, 1]$,
\[
\lim_{n \to \infty}\|\pi^\rho_n-\mu_n\times \mu_n\|_\TV=0,
\]
i.e., the $\mu$-random walk on $\Gamma$ is noise sensitive in total variation.
\end{theorem}

For the definition of affine Weyl groups, see Section \ref{Sec:affineWeyl}.
The infinite dihedral group 
\[
D_\infty:=\abr{s_1, s_2 \mid s_1^2=s_2^2=\id} \quad \text{with $S=\{s_1, s_2\}$}
\]
is an example of affine Weyl group called type $\wt A_1$.
See more examples of affine Weyl groups in Section \ref{Sec:example}.
In fact, there exist a constant $C>0$ and an integer $m>0$ such that for all integer $n>1$,
\[
\|\pi^\rho_n-\mu_n\times \mu_n\|_\TV \le \frac{C(\log n)^m}{\sqrt{n}}.
\]
This strengthening of Theorem \ref{Thm:NSonWeylmain} is stated below as Theorem \ref{Thm:NSonWeyl}.

In Theorem \ref{Thm:NSonWeylmain},
the laziness (i.e., $\mu(\id)>0$) is crucial since otherwise the random walk on that group is not necessarily noise sensitive.
This in particular shows that the noise sensitivity is a property of the random walk rather than the group as was pointed out in \cite{BenjaminiBrieussel}.
Let us note that the $\mu$-random walk on the infinite dihedral group considered there has a particular form: $\mu(s_1)=\mu(s_2)=\mu(\id)=1/3$. 
It is not clear from their proof whether changing the laziness (i.e., the measure on the identity element) would still provide a noise sensitive random walk or not.
We show that this is indeed the case, furthermore, $\mu$ is allowed to be a non-uniform distribution on $S\cup \{\id\}$.


In general,
the following is known for a finitely supported probability measure $\mu$:
If either the group $\Gamma$ admits a surjective homomorphism onto $\Z$, or $(\Gamma, \mu)$ is {\bf non-Liouville}, i.e., there exists a non-constant bounded $\mu$-harmonic function on $\Gamma$,
then a $\mu$-random walk on $\Gamma$ is {\it not} noise sensitive \cite[Theorem 1.1 (2)]{BenjaminiBrieussel}.
In a more specific class of $(\Gamma, \mu)$ (where possibly $\Gamma$ does or does not admit a surjective homomorphism onto $\Z$), 
a strong negation of noise sensitivity has been shown for non-elementary word hyperbolic groups, e.g., free groups of rank at least $2$.
If $\mu$ is non-elementary and has a finite first moment,
then there exists a $\rho_0 \in (0, 1]$ such that
$\|\pi^\rho_n -\mu_n\times\mu_n\|_\TV \to 1$ as $n \to \infty$ for all $\rho \in [0, \rho_0)$
\cite[Theorem 1.3]{NonNS}.
In this generality, it is not known as to whether $\rho_0=1$ or not. 
Note that it is rather straightforward to check that $\rho_0=1$ for simple random walks on free {\em semi}-groups of rank at least $2$, cf.\ \cite[Introduction]{NonNS}.

\subsection*{Outline of the proof of Theorem \ref{Thm:NSonWeylmain}}
An affine Weyl group $(\Gamma, S)$ with the standard set of generators $S$ is associated with some Euclidean space $\R^m$.
The group $\Gamma$ acts on $\R^m$ isometrically and generators act as reflections relative to hyperplanes.
The group $\Gamma$ has the form $\Lambda \rtimes W$ where $\Lambda$ is identified with a lattice in $\R^m$ and $W$ is a finite group (called a spherical Weyl group).
There is a point $o$ in $\R^m$ such that the orbit map $\F: x \mapsto x.o$ is injective.
Taking a conjugate by a translation if necessary, we assume that $o$ is the origin in $\R^m$
(see Section \ref{Sec:affineWeyl} for the precise discussion).
A main ingredient is to establish a local central limit theorem (Theorem \ref{Thm:LCLT}).
We define a discrete normal distribution $\Nc^\F_{n\SS}$ on $\Gamma$ induced from the normalized restriction on $\F(\Gamma)$ of the $m$-dimensional Gaussian density function with a covariance matrix $n\SS$. Further we show the following:
The distribution $\mu_n$ is approximated by $\Nc^\F_{n\SS}$ uniformly on $\Gamma$ within an error of order $n^{-\frac{m+1}{2}}$ as $n$ tends to infinity.
The local central limit theorem itself follows from a classical argument based on characteristic functions.
Some more general results (other than the Cayley graphs of affine Weyl groups) have been proved, e.g., in \cite{KramliSzasz}, \cite{PollicottSharpRates}, \cite{KotaniShiraiSunada} and \cite{SunadaTopological}.
We provide the proof of the local central limit theorem in our setting with an error estimate.
Furthermore, we use an explicit form of matrix $\SS$ in terms of harmonic $1$-forms on a finite quotient graph of the Cayley graph by the lattice $\Lambda$.
The matrix $\SS$ is obtained as some limiting form, which has previously appeared in the literature of symbolic dynamics, see e.g.\ \cite{PollicottSharpRates}; however, the explicit form in Theorem \ref{Thm:LCLT} plays an important role.

We apply this discussion to $(\Gamma \times \Gamma, \pi^\rho)$ for $\rho \in (0, 1]$.
The local central limit theorem enables us to show that there exists a constant $C>0$ such that 
for all large enough $n$, 
\[
\|\pi^\rho_n-\Nc^\F_{n\SS^\rho}\|_\TV \le \frac{C(\log n)^m}{\sqrt{n}}.
\] 
See Theorem \ref{Thm:LCLT_ell1}.
(For the discussion on the sharpness of this bound, see Remark \ref{Rem:sharp}.)
If $\mu$ has support $S\cup \{\id\}$, then $\pi^\rho$ has support $(S\cup \{\id\})^2$ consisting of elements of order at most $2$.
The explicit formula of the covariance matrix implies that $\SS^\rho$ has a block diagonal form and $\SS^\rho=\SS^1$ for all $\rho \in (0, 1]$.
(This is the only part where we use the particular structure of the generating set.)
Thus by the triangle inequality we conclude Theorem \ref{Thm:NSonWeylmain} (in Theorem \ref{Thm:NSonWeyl}).
Let us note that if the support of $\mu$ does not contain $\id$, then the support of $\pi^\rho$ does not generate the group $\Gamma \times \Gamma$ (cf.\ Section \ref{Sec:NSonWeyl}).

\subsection*{Organization}
In Section \ref{Sec:pre},
we introduce affine Weyl groups and discuss background.
In Section \ref{Sec:LCLTsec},
we show the local central limit theorem in a slightly extended setting (Theorem \ref{Thm:LCLT}).
In Section \ref{Sec:applications},
 we deduce noise sensitivity for affine Weyl groups (Theorem \ref{Thm:NSonWeyl}), presenting explicit examples (Section \ref{Sec:example}).
In Appendix \ref{Sec:appendix}, we include the result (Theorem \ref{Thm:Z}) on $\Z^m$ for the sake of convenience.

\subsection*{Notations}
For a constant $C$, we write $C=C_\SS$
to indicate its dependence on $\SS$.
For non-negative real valued functions $f$ and $g$ on a common (sub-)set of non-negative integers $\Z_+$, 
we write $f(n)=O(g(n))$ or $f\ll g$ if there exists a constant $C$ such that $f(n)\le C g(n)$ for all large enough $n$.
We also write $f(n)=O_\SS(g(n))$ if $C=C_{\SS}$ in the above notation.
Further we write $f(n)=\Omega(g(n))$ if there exists a constant $c>0$ such that $f(n) \ge c g(n)$ for all large enough $n$, and $f(n)=\Theta(g(n))$ if $f\ll g$ and $g\ll f$.
For a set $A$, we denote by $|A|$ the cardinality.

\section{Preliminaries}\label{Sec:pre}

For a group $\Gamma$ and for a subset $A$ in $\Gamma$,
we write $\Gamma=\abr{A}$ if $\Gamma$ is generated by $A$ as a {\em semigroup}, i.e., every element in $\Gamma$ is obtained as a product of some finite sequence of elements from $A$.
Let $\Gamma$ be a finitely generated group with a finite symmetric set of generators $S$, i.e., $\Gamma=\abr{S}$ and $S$ is invariant under the map $s \mapsto s^{-1}$.
It holds that, in fact, $s=s^{-1}$ for all $s \in S$ if $S$ consists of involutions.
(This is the case of an (affine) Weyl group $(\Gamma, S)$ in the following discussion.)
Let $\Cay(\Gamma, S)$ be the (right) Cayley graph of $\Gamma$ with respect to $S$, i.e., the set of vertices is $\Gamma$ and an edge $\{x, y\}$ is defined if and only if $x^{-1}y \in S$.
Since $S$ is invariant under $s \mapsto s^{-1}$, 
the Cayley graph is defined as an undirected graph.
For $x \in \Gamma$, let $|x|_S$ denote the word norm with respect to $S$, i.e., the graph distance between $\id$ and $x$ in $\Cay(\Gamma, S)$.

\subsection{Affine Weyl groups}\label{Sec:affineWeyl}

Let $(\Gamma, S)$ be an {\bf affine Weyl group} where $S$ is a canonical finite set of generators,
consisting of involutions, i.e., $s^2=\id$ for every $s \in S$.
The group $\Gamma$ admits a semi-direct product structure $\Gamma=\Lambda \rtimes W$ where
the subgroup $W$ called the {\bf spherical Weyl group} is finite and the normal subgroup $\Lambda$ is isomorphic to a free abelian group of finite rank.
For a thorough background on the subject, we refer to \cite{AbramenkoBrown}.
We also refer to Section \ref{Sec:example} for the examples most relevant to the present discussion.

The group $\Gamma$ is equipped with an isometric action on the standard Euclidean space $\R^m$ for some $m\ge 1$,
where 
each generator $s \in S$ acts as a reflection with respect to an affine hyperplane.
The action is properly discontinuous and admits a relatively compact convex fundamental domain with nonempty interior $C_0$ called a {\bf chamber}.
The group $\Gamma$ acts on the set of chambers $\Cc:=\{x C_0\}_{x \in \Gamma}$ simply transitively, i.e., for all $C_1, C_2 \in \Cc$ there exists $x \in \Gamma$ such that $C_1=xC_2$, and if $x C_0=C_0$, then $x=\id$.
The normal subgroup $\Lambda$ acts freely (i.e., without fixed points) as translations on $\R^m$.
The $\Lambda$-orbit of the origin is a lattice:
\[
\big\{ a_1 v_1+\cdots+a_m v_m : a_1, \dots, a_m \in \Z\big\},
\]
where $v_1, \dots, v_m$ form a basis in $\R^m$.
We identify $\Lambda$ with the lattice in $\R^m$.
Note that the action of $W$ preserves $\Lambda$.

The affine Weyl group $(\Gamma, S)$ is called {\bf reducible} if there exist nontrivial affine Weyl groups generated by  $S_1$ and $S_2$ respectively with $S=S_1\times \cbr{\id} \cup \cbr{\id}\times S_2$ for which $\Gamma=\abr{S_1} \times \abr{S_2}$, and {\bf irreducible} otherwise.
The group $\Gamma$ we consider is possibly (and basically) reducible.
Irreducible ones are completely classified in terms of root systems.
For example, the affine Weyl group of type $\wt A_1$ is the infinite dihedral group
\[
\abr{s_1, s_2 \mid s_1^2=s_2^2=\id},
\]
where $s_1$ and $s_2$ act as reflections with respect to $0$ (the origin) and $1$ respectively in $\R$. 
A chamber has the form of interval $[0, 1)$.

Let us fix a point $o$ in the interior of a chamber and define
\[
\F: \Gamma \to \R^m, \quad x \mapsto x.o.
\]
The map $\F$ is injective since $\Gamma$ acts on the set of chambers simply transitively, and is $\Gamma$-equivariant, i.e., $\F(xy)=x.\F(y)$ for all $x, y \in \Gamma$.
Let us call $\F:\Gamma \to \R^m$ an associated {\bf equivariant embedding}.
Since the generators act as reflections with respect to affine hyperplanes, 
it is illustrative to consider that $\Cay(\Gamma, S)$ is realized in $\R^m$ via the map $\F$.
Namely, the vertices are placed inside of the chambers as the orbit $x.o$ for $x \in \Gamma$ and an edge is a line segment connecting two vertices for which one is obtained from the other by a reflection of the form $x s x^{-1}$ for $s \in S$ and $x \in \Gamma$.

The group $\Lambda$ itself acts on $\Cay(\Gamma, S)$ from left freely as automorphisms of the graph.
Let us consider the quotient graph 
$G=\Lambda\backslash \Cay(\Gamma, S)$.
The graph $G=(V(G), E^{un}(G))$ is finite, the set of vertices $V(G)$ is $W$ and the set of edges $E^{un}(G)$ consists of {\em undirected} edges.
Note, however, that $G$ is {\em not} the right Cayley graph of $W$ with respect to the image $\wbar S$ of $S$ under the quotient map $\Lambda\rtimes W \to W$.
This is because the quotient map restricted on $S$ is not bijective onto $\wbar S$.
The graph $G$ has possibly multiple edges.

\subsection{Pointed finite networks as quotients}\label{Sec:network}

The main interest is on an affine Weyl group $\Gamma$ and the canonical set of generators $S$.
It is, however, useful to discuss a slightly more general setting.
Let $\Gamma$ be a virtually finite rank free abelian group, i.e., $\Gamma$ admits a finite rank free abelian group $\Lambda$ as a finite index subgroup.
We assume that $\Gamma$ acts on an $\R^m$ isometrically with a relatively compact fundamental domain with nonempty interior, and that $\Lambda$ acts as translations and is identified with a lattice in $\R^m$.
Let us fix a point $o$ in the interior of such a fundamental domain of $\Gamma$ and define $\F:\Gamma \to \R^m$ by $x \mapsto x.o$.
The map $\F$ is equivariant with $\Gamma$-actions and injective.
For a finite symmetric set of generators $S$ in $\Gamma$,
let 
\[
G:=\Lambda \backslash \Cay(\Gamma, S).
\]
The quotient $G=(V(G), E^{un}(G))$ is a finite (undirected) graph possibly with multiple edges (whence a multi-graph) and with loops.
It holds that $V(G)=\Lambda \backslash \Gamma$
and 
\[
E^{un}(G)=\{\{x, x.s\} \ : \ x \in \Lambda \backslash \Gamma, \ s \in S\},
\]
where $\{x, x.s\}$ and $\{x.s, x\}$ are identified.
For the examples, see Section \ref{Sec:example}.

Let $\mu$ be a probability measure on $\Gamma$ such that the support $\supp \mu$ of $\mu$ is finite, that $\Gamma=\abr{\supp \mu}$, and that $\mu$ is {\bf symmetric}, i.e., $\mu(s)=\mu(s^{-1})$ for every $s \in \supp \mu$.
If we define $S=\supp \mu$, then $S$ is a finite symmetric set of generators.
Let 
\[
p(\{x, x.s\}):=\mu(s) \quad \text{for $\{x, x.s\} \in E^{un}(G)$}.
\]
This defines a Markov chain on $G=\Lambda \backslash \Cay(\Gamma, S)$ with transition probabilities 
\[
\sum_{s \in S, y=x.s}\mu(s) \quad \text{for $x, y \in \Lambda\backslash \Gamma$}.
\]
Note that this Markov chain is irreducible, i.e., it visits every vertex from every other vertex after some time since $\Gamma=\abr{\supp \mu}$.
Furthermore, it
is reversible with respect to the uniform distribution $\pi$ on the set of vertices $V(G)=\Lambda \backslash \Gamma$.
Indeed, since $\mu$ is symmetric,
it holds that
\begin{equation}\label{Eq:reversible}
\pi(x) p(\{x, x.s\})=\pi(x.s) p(\{x.s, x\}) \quad \text{for $\{x, x.s\} \in E^{un}(G)$},
\end{equation}
where $\pi(x)=1/|V(G)|$ for $x \in V(G)$.
Note that if $\mu(\id)>0$, then $p(x, x) \ge \mu(\id)>0$ for every $x \in V(G)$.
For each $\{x, y\} \in E^{un}(G)$, let us define the {\bf conductance} by
\[
c(\{x, y\}):=\pi(x)p(\{x, y\}).
\]
This is well-defined since $c(\{x, y\})=c(\{y, x\})$ by \eqref{Eq:reversible}.
Note that $c(\{x, y\})>0$ for all $\{x, y\} \in E^{un}(G)$.
Let $x_0 \in \Lambda \backslash \Gamma$ denote the coset containing $\id$. 
Let us call $(G, c, x_0)$ the pointed (finite) {\bf network} as the finite multi-graph $G$ equipped with the conductance $c:E^{un}(G) \to (0, \infty)$ and the point $x_0$.

\section{Local central limit theorems}\label{Sec:LCLTsec}

\subsection{Harmonic $1$-forms on finite graphs}\label{Sec:harmonic}

Let $(G, c, x_0)$ be the pointed finite network.
Henceforth it is convenient to consider $G$ as a {\bf graph with orientations} where 
each edge (and loop) has both possible orientations.
Let
\[
E(G):=\big\{(x, y), (y, x) \ : \ \{x, y\} \in E^{un}(G)\big\}.
\]
For $e=(x, y)$, we write $\wbar e=(y, x)$.
The ``reversing direction'' operation $\wbar{\ \cdot\ }: E(G) \to E(G)$, $e\mapsto \wbar e$, defines a bijection and $\wbar{\wbar e}=e$ for $e \in E(G)$.
For $e=(x, y) \in E(G)$, let us denote by $oe:=x$ the {\bf origin} and by $te:=y$ the {\bf terminus} of $e$ respectively.
We have that $o\wbar e=te$ for $e \in E(G)$.
Let us also consider $\Cay(\Gamma, S)$ as a graph with orientations, defining both possible orientations for each edge and loop.
Letting $c(e):=c(x, y)$ and $p(e):=p(x, y)$ for $e=(x, y)$,
we have that
$c(e)=c(\wbar e)$ and $c(e)=\pi(oe)p(e)$ for $e \in E(G)$.
It holds that by the definition of conductance, 
\[
\pi(x)=\sum_{e \in E_x} c(e), \quad \text{where $E_x:=\big\{e \in E(G) \ : \ oe=x\big\}$ for $x \in V(G)$}.
\]

Let us define the $\C$-linear space of complex-valued functions on $V(G)$ by
\[
C^0(G, \C):=\big\{f: V(G) \to \C\big\}
\]
equipped with the inner product
$\abr{f_1, f_2}_\pi:=\sum_{x \in V(G)}f_1(x)\wbar{f_2(x)}\pi(x)$,
where $\wbar{a}$ stands for the complex-conjugate of $a \in \C$.
Similarly, let $C^0(G, \R)$ be the $\R$-linear space of real-valued functions on $V(G)$ endowed with the inner product as the restriction of $\abr{\cdot, \cdot}_\pi$.
Further let us define the $\R$-linear space of real-valued {\bf $1$-forms} on $E(G)$ by
\[
C^1(G, \R):=\big\{\omega: E(G) \to \R \ : \ \text{$\omega(\wbar e)=-\omega(e)$ for $e \in E(G)$}\big\}
\]
equipped with the inner product
$\abr{\omega_1, \omega_2}_c:=(1/2)\sum_{e \in E(G)}\omega_1(e)\omega_2(e)c(e)$.
The {\bf differential} $d: C^0(G, \R) \to C^1(G, \R)$ is the $\R$-linear map defined by
\[
df(e):=f(te)-f(oe) \quad \text{for $e \in E(G)$}.
\]
Moreover, the {\bf adjoint} $d^\ast:C^1(G, \R) \to C^0(G, \R)$ with respect to the inner products is obtained by
\[
d^\ast \omega(x):=-\sum_{e \in E_x}\frac{1}{\pi(x)}c(e) \omega(e) \quad \text{for $x \in V(G)$}.
\]
It holds that for $f \in C^0(G, \R)$ and $\omega\in C^1(G, \R)$,
\begin{equation}\label{Eq:adjoint}
\abr{df, \omega}_c=\abr{f, d^\ast \omega}_\pi.
\end{equation}
Note that if we define the transition operator $P$ on $C^0(G, \R)$ to itself by
\[
Pf(x):=\sum_{e \in E_x}\frac{1}{\pi(x)}c(e)f(te) \quad \text{for $x \in V(G)$},
\]
then 
$d^\ast d=I-P$ where $I$ is the identity operator.
Let us define the space of {\bf harmonic $1$-forms} by
\[
H^1:=\big\{\omega \in C^1(G, \R) \ : \ d^\ast \omega=0\big\}.
\]
Note that $H^1=(\Im d)^\bot$ the orthogonal complement of the image $\Im d$ by \eqref{Eq:adjoint}.
The fact that we use in the sequel is that for every $1$-form $\omega \in C^1(G, \R)$ there exists a {\em unique} harmonic $1$-form $u \in H^1$ and some $f \in C^0(G, \R)$ such that
\[
u+df=\omega.
\]
The $u$ is obtained as the $H^1$-part in the orthogonal decomposition $C^0(G, \R)=H^1\oplus \Im d$.
Note that $f$ is not unique since every $f$ added a constant function satisfies the relation.

\begin{remark}
If we endow $G$ with a structure of $1$-dimensional CW complex, then the $1$-cohomology group with real coefficient is defined as $H^1(G, \R):=C^1(G, \R)/\Im d$.
The fact mentioned above means that every $1$-cohomology class is represented by a unique harmonic $1$-form.
Although all these notions are not needed in our discussion, it might be useful to grasp an idea behind some of our computations.
\end{remark}

\subsection{Perturbations of transfer operators}\label{Sec:transfer}

For each $1$-form $\omega \in C^1(G, \R)$,
let us define the {\bf transfer operator} on $C^0(G, \C)$ by
\[
\Lc_\omega f(x):=\sum_{e \in E_x}p(e) e^{2\pi i \omega(e)}f(te) \quad \text{for $x \in V(G)$}.
\]
Here $i=\sqrt{-1}$.
In this particular setting where $\mu$ is symmetric,
the transfer operator is self-adjoint on $(C^0(G, \C), \abr{\cdot, \cdot}_\pi)$.
Hence it has real eigenvalues.
Let $\lambda(\omega)$ be the largest eigenvalue of $\Lc_\omega$.
If $\omega=0$, then $\lambda(0)=1$ and this is a simple eigenvalue by the Perron-Frobenius theorem since $\Lc_0=P$ and $P$ is irreducible.
We apply to an analytic perturbation in $\omega$:
for a small enough neighborhood $U$ of $0$ in $C^1(G, \R)$,  
the function $U \to \R$, $\omega \mapsto \lambda(\omega)$ is real analytic.
Moreover, corresponding eigenvectors $f_\omega$ depend analytically in $\omega \in U$ with $f_0=\1$ (the constant vector with all $1$'s).
This follows from the implicit function theorem for $\det(tI-\Lc_\omega)=0$ around $(t, \omega)=(1, 0)$ in this finite graph setting.
Note that if we consider $\omega+d\f$ for $\f \in C^0(G, \R)$ in place of $\omega$,
then $\Lc_{\omega+d\f}=e^{-2\pi i\f}\Lc_\omega e^{2\pi i\f}$, where $(e^\f f)(x):=e^{\f(x)}f(x)$ for $x \in V(G)$, and thus
\[
\lambda(\omega+d\f)=\lambda(\omega).
\]
This shows that $\lambda(\omega)$ depends only on the harmonic $H^1$-part of $\omega$.
Let
\[
\beta(\omega):=\log \lambda(\omega).
\]
Since $\Lc_0=P$ which has a simple eigenvalue $\lambda(0)=1$, there exists
a small enough open neighborhood $U$ of $0$ in $C^1(G, \R)$ such that $\lambda(\omega)$ is a simple eigenvalue of $\Lc_{\omega}$ and $\beta(\omega)=\log \lambda(\omega)$ is well-defined for all $\omega \in U$.
It holds that $\beta(0)=0$ since $\lambda(0)=1$.

\begin{lemma}\label{Lem:der}
Let $\lambda(\omega)=e^{\beta(\omega)}$ and $f_\omega$ be the eigenvalue and the corresponding eigenvector of $\Lc_\omega$ such that $f_0=\1$ and $\abr{f_\omega, \1}_\pi=1$ for $\omega \in U$ where $U$ is a neighborhood of $0$ in $C^1(G, \R)$.
For all harmonic $1$-forms $u, u_i \in H^1$ on $G$ and real parameters $r, r_i$ for $i=1,2,3$, the following holds:
\begin{equation}\label{Eq:Lem:der1}
\frac{d}{dr}\Big|_{r=0}\beta(ru)=0,
\end{equation}
\begin{equation}\label{Eq:Lem:eigenvector}
\frac{d}{dr}\Big|_{r=0}f_{ru}(x)=0 \quad \text{for all $x \in V(G)$},
\end{equation}
\begin{equation}\label{Eq:Lem:der2}
\frac{\partial^2}{\partial r_1 \partial r_2}\Big|_{(r_1, r_2)=(0, 0)}\beta(r_1 u_1+r_2 u_2)=-4\pi^2 \sum_{e \in E(G)}u_1(e)u_2(e)c(e),
\end{equation}
and
\begin{equation}\label{Eq:Lem:der3}
\frac{\partial^3}{\partial r_1 \partial r_2 \partial r_3}\Big|_{(r_1, r_2, r_3)=(0, 0, 0)}\beta(r_1u_1+r_2u_2+r_3u_3)=0.
\end{equation}
\end{lemma}

\proof
Let $f_\omega$ be the eigenvector normalized as stated.
Since $\wbar{\Lc_\omega f_\omega}=\Lc_{-\omega}\wbar{f_\omega}$ holds for all $\omega$ and $\lambda(\omega)$ is real for all $\omega$ in a small enough neighborhood of $0$ in $C^1(G, \R)$,
it holds that $\lambda(-\omega)=\wbar{\lambda(\omega)}=\lambda(\omega)$.
Thus $\beta(-\omega)=\beta(\omega)$ for all $\omega \in U$ and all odd time derivatives of $\beta$ at $0$ vanish.
This in particular implies \eqref{Eq:Lem:der1} and \eqref{Eq:Lem:der3}.

For every $1$-form $u$, it holds that $ru \in U$ for all small enough real $r$ and
\begin{equation}\label{Eq:normalization_der}
\left\langle\frac{d}{d r}\Big|_{r=0}f_{ru}, \1\right\rangle_\pi=0.
\end{equation}
Moreover, since $P$ is self-adjoint with respect to $\abr{\cdot,\cdot}_\pi$ and $P\1=\1$, by \eqref{Eq:normalization_der} it holds that
\begin{equation}\label{Eq:normalization_der_P}
\left\langle P\(\frac{d}{d r}\Big|_{r=0}f_{ru}\), \1\right\rangle_\pi=\left\langle \frac{d}{d r}\Big|_{r=0}f_{ru}, P\1\right\rangle_\pi=\left\langle\frac{d}{d r}\Big|_{r=0}f_{ru}, \1\right\rangle_\pi=0.
\end{equation}
We will also use the analogous identities to \eqref{Eq:normalization_der} and \eqref{Eq:normalization_der_P} for the second derivatives of the normalization: $\abr{f_\omega, \1}_\pi=1$ for $\omega \in U$.

First
differentiating $\Lc_{r u}f_{ru}=e^{\beta(ru)}f_{ru}$ at $r=0$ yields for each $x \in V(G)$,
\begin{equation}\label{Eq:Lem:point}
\sum_{e \in E_x}\(p(e)(2\pi i u(e))+p(e)\frac{d}{d r}\Big|_{r=0}f_{r u}(te)\)=
\frac{d}{d r}\Big|_{r=0}\beta(ru)+\frac{d}{d r}\Big|_{r=0}f_{ru}(x),
\end{equation}
where we have used $f_0=\1$.
Let us note that \eqref{Eq:Lem:point} above yields by \eqref{Eq:Lem:der1} which we have just shown and by that $d^\ast u=0$,
\[
P\(\frac{d}{d r}\Big|_{r=0}f_{ru}\)(x)=\frac{d}{d r}\Big|_{r=0}f_{ru}(x)
\quad \text{for each $x \in V(G)$}.
\]
Since $P$ has the simple eigenvalue $1$, this implies that $(d/d r)|_{r=0}f_{ru}$ is constant.
By \eqref{Eq:normalization_der}, for every harmonic $1$-from $u$,
it holds that 
$(d/dr)|_{r=0}f_{ru}(x)=0$ for all $x \in V(G)$,
showing \eqref{Eq:Lem:eigenvector}.

For all $1$-forms $u_1$ and $u_2$ and for all small enough reals $r_1$ and $r_2$,
it holds that
\[
\Lc_{r_1 u_1+r_2 u_2}f_{r_1, r_2}=e^{\beta_{r_1, r_2}}f_{r_1, r_2}
\quad \text{where $\beta_{r_1, r_2}:=\beta(r_1u_1+r_2u_2)$ and $f_{r_1, r_2}:=f_{r_1 u_1+r_2 u_2}$}.
\]
Taking the second derivatives at $(r_1, r_2)=(0, 0)$ of both terms yields by \eqref{Eq:Lem:der1}, \eqref{Eq:Lem:eigenvector} and that $f_0=\1$, for each $x \in V(G)$,
\begin{align*}
&\sum_{e \in E_x}\Bigg(p(e)(-4\pi^2 u_1(e)u_2(e))+p(e)\frac{\partial^2}{\partial r_1 \partial r_2}\Big|_{(r_1, r_2)=(0, 0)}f_{r_1, r_2}(te)\Bigg)\\
&\qquad \qquad \qquad \qquad \qquad \qquad=\frac{\partial^2}{\partial r_1 \partial r_2}\Big|_{(r_1, r_2)=(0, 0)}\beta_{r_1, r_2}
+\frac{\partial^2}{\partial r_1 \partial r_2}\Big|_{(r_1, r_2)=(0, 0)}f_{r_1, r_2}(x).
\end{align*}
Evaluating the inner products of the above terms with $\1$ leads
\begin{align*}
&-4\pi^2\sum_{e \in E(G)}c(e)u_1(e)u_2(e)
+\left\langle P\(\frac{\partial^2}{\partial r_1 \partial r_2}\Big|_{(r_1, r_2)=(0, 0)}f_{r_1, r_2}\), \1\right\rangle_\pi\\
&\qquad \qquad \qquad \qquad \qquad \qquad=\frac{\partial^2}{\partial r_1 \partial r_2}\Big|_{(r_1, r_2)=(0, 0)}\beta_{r_1, r_2}
+\left\langle \frac{\partial^2}{\partial r_1 \partial r_2}\Big|_{(r_1, r_2)=(0, 0)}f_{r_1, r_2}, \1\right\rangle_\pi.
\end{align*}
The second terms in the left hand side and in the right hand side respectively are $0$
since they are the second derivatives on the normalization (cf.\ \eqref{Eq:normalization_der} and \eqref{Eq:normalization_der_P}). 
This proves \eqref{Eq:Lem:der2}.
\qed

\subsection{An explicit Hessian formula in terms of harmonic $1$-forms}\label{Sec:Hessian}

Recall that $\F:\Gamma \to \R^m$, $x \mapsto x.o$.
Taking a conjugate to the action of $\Gamma$ by a translation in $\R^m$, we assume that $o$ is the origin, whence $\F(\id)=0$.
The function $(x, y) \mapsto \F(y)-\F(x)$ for (oriented) edges $(x, y)$ in $\Cay(\Gamma, S)$ is invariant under the action by $\Lambda$.
Indeed, this follows since $\Lambda$ is identified with a lattice and acts as translations in $\R^m$.
Therefore this descends to an $\R^m$-valued function on $E(G)$, which we denote by 
$e \mapsto \F_e$ for $e \in E(G)$.
Note that $\F_{\wbar e}=-\F_e$ for each $e \in E(G)$.
Let $\abr{\cdot, \cdot}$ be the standard inner product in $\R^m$.
For each $v \in \R^m$, 
let
\[
\wh v(e):=\abr{v, \F_e} \quad \text{for $e \in E(G)$}.
\]
This $\wh v$ defines a $1$-form on $G$.
For $v \in \R^m$ near $0$,
let $\beta(v):=\beta(\wh v)$ where $e^{\beta(v)}$ is the largest eigenvalue of the transfer operator $\Lc_{\wh v}$.
Let us define the Hessian of $\beta$ at $0$ in $\R^m$ with the standard coordinate $(r_1, \dots, r_m)$ by
\[
\Hess_0\beta:=\(\frac{\partial^2}{\partial r_k\partial r_l}\Big|_{(r_1, \dots, r_m)=(0, \dots, 0)}\beta(r_1, \dots, r_m)\)_{k, l=1, \dots, m}.
\]
For the pointed finite network $(G, c, x_0)$ and $\F:\Gamma \to \R^m$, we compute $\Hess_0\beta$.
\begin{lemma}\label{Lem:H}
The Hessian $\Hess_0 \beta$ of $\beta$ at $0$ in $\R^m$ is non-degenerate and negative definite.
Moreover, it holds that
\begin{equation}\label{Eq:vSSv}
\abr{v_1, \Hess_0 \beta\, v_2}=-4\pi^2\sum_{e \in E(G)}u_1(e)u_2(e)c(e),
\end{equation}
where $u_i$ is defined as the harmonic part of $\wh v_i$ for $v_i \in \R^m$, $i=1, 2$. 
\end{lemma}

\proof
For every $v \in \R^m$,
we have $\wh v(e)=\abr{v, \F_e}$ for $e \in E(G)$,
and $u$ is the harmonic part of $\wh v$, i.e., the unique $u \in H^1$ such that
$u+df=\wh v$ for some $f \in C^0(G, \R)$.
Note that the resulting map $v \mapsto u$ is $\R$-linear.
Since $\beta(v)$ depends only on the harmonic part of $v$, Lemma \ref{Lem:der} \eqref{Eq:Lem:der2} implies that
\begin{align*}
\abr{v_1, \Hess_0\beta\, v_2}
&=\frac{\partial^2}{\partial r_1 \partial r_2}\Big|_{(r_1, r_2)=(0, 0)}\beta(r_1 v_1+r_2 v_2)\\
&=\frac{\partial^2}{\partial r_1 \partial r_2}\Big|_{(r_1, r_2)=(0, 0)}\beta(r_1 u_1+r_2 u_2)=-4\pi^2\sum_{e \in E(G)}u_1(e)u_2(e) c(e).
\end{align*}
This shows \eqref{Eq:vSSv}.

Let $v_1, \dots, v_m$ be a basis of the lattice in $\R^m$:
$\Lambda=\big\{ a_1 v_1+\cdots+a_m v_m : a_1, \dots, a_m \in \Z\big\}$.
For each $v_k=\F(v_k) \in \Lambda$ under the identification of $\Lambda$ with the lattice, 
there exists a path $(\wt e_1, \dots, \wt e_n)$ from $\id$ to $v_k$ in $\Cay(\Gamma, S)$
since the Cayley graph is connected.
Let $(e_1, \dots, e_n)$ be the image in $G$ of that path under the covering map from $\Cay(\Gamma, S)$.
Note that the image is a cycle: $x_0=oe_1$, $te_i=oe_{i+1}$ for $i=1, \dots, n-1$ and $te_n=x_0$.
Thus $\sum_{l=1}^n df(e_l)=0$ and
$\sum_{l=1}^n u(e_l)=\sum_{l=1}^n (u(e_l)+df(e_l))=\sum_{l=1}^n \wh v(e_l)$.
Furthermore,
\begin{align*}
\sum_{l=1}^n \wh v(e_l)=\sum_{l=1}^n \abr{v, \F_{e_l}}
=\sum_{l=1}^n \abr{v, \F(t\wt e_l)-\F(o\wt e_l)}
=\abr{v, \sum_{l=1}^n(\F(t\wt e_l)-\F(o\wt e_l))}.
\end{align*}
The last term equals $\abr{v, v_k}$
since $\sum_{l=1}^n (\F(t\wt e_l)-\F(o\wt e_l))=\F(t\wt e_n)-\F(o \wt e_1)=v_k$.
This shows the following: 
For each $k=1, \dots, m$ there exists a cycle $(e_1, \dots, e_n)$ in $G$ with $x_0=oe_1$ and $te_n=x_0$ such that
\begin{equation}\label{Eq:period}
\sum_{l=1}^n u(e_l)=\abr{v, v_k}.
\end{equation}

For $v \in \R^m$,
let us assume that $\abr{v, \Hess_0\beta\, v}=0$.
It holds that $u=0$ by \eqref{Eq:vSSv}, and thus $\abr{v, v_k}=0$ for every $k=1, \dots, m$ by \eqref{Eq:period}. Hence $v=0$ since $v_1, \dots, v_m$ form a basis of a lattice in $\R^m$.
This shows that $\Hess_0\beta$ is non-degenerate. 
Furthermore $\Hess_0\beta$ is negative definite by \eqref{Eq:vSSv}.
\qed

\begin{remark}\label{Rem:cohomology}
Let us consider the Hessian $\Hess_{H^1}\beta$ of $\beta$ at $0$ on $H^1$, i.e.,
\[
\abr{u_1, \Hess_{H^1}\beta\, u_2}_\pi=\frac{\partial^2}{\partial r_1 \partial r_2}\Big|_{(r_1, r_2)=(0, 0)}\beta(r_1 u_1+r_2 u_2) \quad \text{for $u_1, u_2 \in H^1$},
\]
where $H^1 \to H^1: u \mapsto \Hess_{H^1}\beta\, u$ defines an $\R$-linear map.
Lemma \ref{Lem:der} \eqref{Eq:Lem:der2} implies that
\[
\abr{u_1, \Hess_{H^1}\beta\, u_2}_\pi=-4\pi^2 \sum_{e \in E(G)}u_1(e)u_2(e)c(e),
\]
which shows that $\Hess_{H^1}\beta$ is non-degenerate and negative definite on $H^1$.
There exists a natural inclusion $\R^m=H^1(\R^m\slash \Lambda, \R) \to H^1(G, \R)$, represented by $\R^m\to H^1: v \mapsto u$ in Lemma \ref{Lem:H}.
In this identification, $\Hess_0\beta$ is the restriction of $\Hess_{H^1}\beta$ to $\R^m$,
and this implies that $\Hess_0\beta$ is non-degenerate and negative definite.
A thorough framework is found in \cite{SunadaTopological}.
We use the explicit form of $\Hess_0\beta$ later in our discussion.
\end{remark}

\subsection{Local central limit theorems}\label{Sec:LCLT}

For every positive integer $n \in \Z_{>0}$, it holds that
\begin{equation}\label{Eq:L}
\Lc_{\wh v}^n \1(x_0)=\sum_{(e_1, \dots, e_n)}p(e_1)\cdots p(e_n)e^{2\pi i (\wh v(e_1)+\cdots+\wh v(e_n))}.
\end{equation}
In the above the summation runs over all directed paths $(e_1, \dots, e_n)$ starting from $x_0$ in $G$, i.e.,
$oe_1=x_0$ and $t e_k=o e_{k+1}$ for each $k=1, \dots, n-1$. 
For each such path $(e_1, \dots, e_n)$, 
there exists a unique path $(\wt e_1, \dots, \wt e_n)$ which is a lift of the path, starting from $\id$ in $\Cay(\Gamma, S)$.
By a lift we mean that the path $(e_1, \dots, e_n)$ is the image of $(\wt e_1, \dots, \wt e_n)$ under the covering map $\Cay(\Gamma, S) \to G=\Lambda\backslash \Cay(\Gamma, S)$.
The definition of the $1$-form $\wh v$ on $G$ implies the following:
\begin{align*}
\wh v(e_1)+\cdots+\wh v(e_n)
&=\abr{v, \F_{e_1}}+\cdots+\abr{v, \F_{e_n}}\\
&=\abr{v, \F(t\wt e_1)-\F(o \wt e_1)}+\cdots+\abr{v, \F(t\wt e_n)-\F(o\wt e_n)}
=\abr{v, \F(t\wt e_n)},
\end{align*}
where $\F(o\wt e_1)=\F(\id)=0$.
Recall that $p(e)=\mu(s)$ for $e=(x, x.s) \in E(G)$ and $s \in S$ in the pointed finite network $(G, c, x_0)$.
By \eqref{Eq:L}, it holds that
\[
\Lc^n_{\wh v}\1(x_0)=\sum_{(\wt e_1, \dots, \wt e_n)}\mu(s_1)\cdots \mu(s_n)e^{2\pi i \abr{v, \F(t\wt e_n)}}
=\sum_{x \in \Gamma}\mu_n(x)e^{2\pi i \abr{v, \F(x)}}.
\]
In the above, the (edge) path $(\wt e_1, \dots, \wt e_n)$ is represented as the (vertex) path on $\id$, $s_1$, $s_1 s_2$, $\dots$, $s_1\cdots s_n$ in $\Cay(\Gamma, S)$.
Therefore letting
\[
\f_{\mu_n}(v):=\sum_{x \in \Gamma}\mu_n(x)e^{2\pi i\abr{v, \F(x)}} \quad \text{for $v \in \R^m$},
\]
we have the following: 
For all $v \in \R^m$ and all $n\in \Z_{>0}$,
\begin{equation}\label{Eq:transfer-char}
\Lc_{\wh v}^n\1(x_0)=\f_{\mu_n}(v).
\end{equation}

Let $\Lambda^\ast$ be the dual lattice of $\Lambda$, i.e., 
\[
\Lambda^\ast:=\Big\{a_1 v_1^\ast+\cdots+a_m v_m^\ast \ : \ a_1, \dots, a_m \in \Z\Big\},
\]
where $v_1^\ast, \dots, v_m^\ast$ form the dual basis of $v_1, \dots, v_m$ in $\R^m$: $\abr{v_k^\ast, v_l}=1$ if $k=l$ and $0$ else.
The fundamental parallelotope of $\Lambda^\ast$ in $\R^m$ is denoted by
\[
D:=\Big\{r_1 v_1^\ast+\cdots+r_m v_m^\ast \in \R^m \ : \ |r_i| \le 1/2, \ i=1, \dots, m\Big\}.
\]
The volume of $D$ is assumed to be $1$ up to a homothety in $\R^m$.
The Fourier inversion formula shows that for all $n \in \Z_{>0}$,
\begin{equation}\label{Eq:Fourier-inversion}
\mu_n(x)=\int_D \f_{\mu_n}(v)e^{-2\pi i \abr{v, \F(x)}}\,dv \quad \text{for $x \in \Gamma$}.
\end{equation}
For $\delta>0$,
let
\[
D_\delta:=\Big\{r_1 v_1^\ast+\cdots+r_m v_m^\ast \in \R^m \ : \ |r_i| < \delta, \ i=1, \dots, m\Big\}.
\]

\begin{lemma}\label{Lem:aperiodic}
If $\mu(\id)>0$, then
for all small enough $\delta>0$,
there exists a constant $c_\delta>0$ such that for all $n \in \Z_{>0}$,
\[
|\f_{\mu_n}(v)|\le \sqrt{|V(G)|}\cdot e^{-c_\delta n} \quad \text{for all $v \in D \setminus D_\delta$}.
\] 
\end{lemma}

\proof
This uses a standard perturbation argument; we provide a proof for the sake of completeness.
For $v \in \R^m$, let $\|\Lc_{\wh v}\|:=\max_{\|f\|_\pi=1}\|\Lc_{\wh v}f\|_\pi$ where $\|\cdot\|_\pi$ is the associated norm in $C^0(G, \C)$. Since $\Lc_{\wh v}$ is self-adjoint, the eigenvalues are real and the operator norm $\|\Lc_{\wh v}\|$ is the spectral radius $|\lambda(v)|$, i.e., the largest eigenvalue in absolute value.
Note that $|\lambda(v)|\le 1$.
For $v \in D$, if $|\lambda(v)|=1$, then the condition $\mu(\id)>0$ implies that $v=0$ (in which case in fact $\lambda(0)\neq -1$).
Indeed, for the maximal eigenvalue $\lambda(v)$ in absolute value and a corresponding eigenvector $f_v$,
we have $\Lc_{\wh v}f_v=\lambda(v)f_v$.
Taking absolute values shows that $|f_v|\le P|f_v|$, implying that $|f_v|$ is a non-zero constant since $P$ is irreducible.
Further since $f_v$ is an eigenvector with the eigenvalue $1$ in absolute value, 
$\lambda(v)f_v(x)$ and $e^{2\pi i \abr{v, \F_e}}f_v(te)$ for $e \in E_x$ are on a common circle in the complex plane for each $x \in V(G)$.
Since $\Lc_{\wh v}f_v=\lambda(v)f_v$,
it holds that for all $x \in V(G)$ and for all $e \in E_x$,
\begin{equation}\label{Eq:circle}
\lambda(v) f_v(x)=e^{2\pi i \abr{v, \F_e}}f_v(te).
\end{equation}
If $\mu(\id)>0$, then for each $v_k \in \Lambda$ there exists an edge path from $\id$ to $v_k$ in $\Cay(\Gamma, S)$ of length with a given (in particular, even) parity.
Applying to \eqref{Eq:circle} along the image of the path in $G$ successively yields $\abr{v, v_k} \in \Z$.
This holds for a basis $v_1, \dots, v_k$ of $\Lambda$, implying that
$v \in \Lambda^\ast$.
Therefore if $v \in D$, then $v=0$.

We have shown that $|\lambda(v)|<1$ for all $v \in D\setminus\{0\}$, and
in this finite dimensional setting, $v \mapsto \|\Lc_{\wh v}\|=|\lambda(v)|$ is continuous.
Thus for a small enough $\delta>0$, there exists a constant $c_\delta>0$ such that $|\lambda(v)| \le e^{-c_{\delta}}$ on a compact set $D \setminus D_\delta$.
Since
\[
|\Lc_{\wh v}^n\1(x_0)| \sqrt{\pi(x_0)} \le \|\Lc_{\wh v}^n\1\|_\pi \le |\lambda(v)|^{n} \|\1\|_\pi,
\]
$\|\1\|_\pi=1$ and $\pi(x_0)=1/|V(G)|$,
by \eqref{Eq:transfer-char}, we conclude the claim.
\qed

\medskip

For an associated $\Gamma$-equivariant embedding $\F:\Gamma \to \R^m$, $x \mapsto x.o$ and a non-degenerate positive definite matrix $\SS$,
let
\[
\xi_{\SS}(x):=\frac{1}{(2\pi)^{\frac{m}{2}}\sqrt{\det \SS}}e^{-\frac{1}{2}\abr{\F(x), \SS^{-1}(\F(x))}} \quad \text{for $x \in \Gamma$}.
\]

\begin{theorem}[Local central limit theorem]\label{Thm:LCLT}
Let $\Gamma$ be a virtually finite rank free abelian group acting on $\R^m$ isometrically with a relatively compact fundamental domain which contains the origin in the interior.
Let $\mu$ be a probability measure on $\Gamma$ such that the support $\supp \mu$ is finite, $\Gamma=\abr{\supp \mu}$ and $\mu$ is symmetric.
If $\mu(\id)>0$,
then the following holds:
There exist a non-degenerate positive definite matrix $\SS$ and a constant $C>0$ 
such that
\[
\sup_{x \in \Gamma}|\mu_n(x)-\xi_{n\SS}(x)|\le \frac{C}{n^{\frac{m+1}{2}}} \quad \text{
for all $n \in \Z_{>0}$}.
\]
Moreover, the matrix $\SS$ is obtained by
\begin{equation}\label{Eq:cov}
\abr{v_1, \SS v_2}=\sum_{e \in E(G)}u_1(e)u_2(e)c(e),
\end{equation}
where $u_i$ is the harmonic part of $\wh v_i$ for $v_i \in \R^m$ for $i=1, 2$. 
\end{theorem}

\proof
For all $\delta>0$, the Fourier inversion formula \eqref{Eq:Fourier-inversion} and the change of variables $v \mapsto v/\sqrt{n}$ yield the following:
For $n \in \Z_{>0}$ and for $x\in \Gamma$,
\begin{align*}
\mu_n(x)&=\int_{D_{\delta}}\f_{\mu_n}\(v\)e^{-2\pi i \abr{v, \F(x)}}\,dv+\int_{D\setminus D_\delta} \f_{\mu_n}(v)e^{-2\pi i \abr{v, \F(x)}}\,dv\\
&=\frac{1}{n^{\frac{m}{2}}}\int_{D_{\delta \sqrt{n}}}\f_{\mu_n}\(\frac{v}{\sqrt{n}}\)e^{-2\pi i \abr{v, \F(x)}/\sqrt{n}}\,dv+\int_{D\setminus D_\delta} \f_{\mu_n}(v)e^{-2\pi i \abr{v, \F(x)}}\,dv.
\end{align*}
Since $\mu(\id)>0$, Lemma \ref{Lem:aperiodic} shows that 
for all small enough $\delta>0$,
there exists a constant $c_\delta>0$ such that for all $n \in \Z_{>0}$,
\[
\sup_{v \in D\setminus D_\delta}|\f_{\mu_n}(v)|=O\(e^{-c_\delta n}\).
\]

We will analyze the first integral in the last displayed equation.
For $v\in \R^m$,
let $\wh v=u+d\f$ and $u$ be the harmonic $H^1$-part of $\wh v$ for some $\f\in C^0(G, \R)$.
There exists a constant $C$ such that such $\f$ can be chosen to satisfy 
$\|\f\|_\infty \le C\|v\|$.
Indeed, letting $v=\sum_{i=1}^m \alpha_i v_i$ for the standard basis $v_1, \dots, v_m$, and $\alpha_1, \dots, \alpha_m \in \R$,
we choose $\f_i$ such that $\wh v_i=u_i+d\f_i$, and define $\f:=\sum_{i=1}^m \alpha_i \f_i$.
This yields the inequality with $C=\sqrt{m}\max_{i=1, \dots, m}\|\f_i\|_\infty$.

We have $\Lc_{\wh v}=e^{-2\pi i\f}\Lc_u e^{2\pi i\f}$.
By \eqref{Eq:transfer-char}, up to replacing $c_\delta$ by a smaller positive value,
\[
\f_{\mu_n}(v)=e^{n\beta(v)}\(\abr{e^{2\pi \f}\1, f_{u}}_\pi e^{-2\pi i\f(x_0)}f_{u}(x_0)+O(e^{-c_\delta n})\) \quad \text{for $v \in D_\delta$},
\]
where $f_u$ is the normalized eigenvector of $\Lc_u$ of eigenvalue $e^{\beta(v)}$.
By Lemma \ref{Lem:H}, the Hessian of $\beta$ at $0$ on $\R^m$ is obtained by 
$\Hess_0\beta=-4\pi^2\SS$.
Lemma \ref{Lem:der} implies that by the Taylor theorem,
\[
\beta(v)=-2\pi^2\abr{v, \SS v}+O(\|v\|^4) \quad \text{for $v \in D_\delta$}.
\]
Therefore for all $n \in \Z_{>0}$ and for all $v \in D_{\delta \sqrt{n}}$,
\[
\beta\(\frac{v}{\sqrt{n}}\)=-\frac{2\pi^2}{n}\abr{v, \SS v}+O\(\frac{\|v\|^4}{n^2}\).
\]
Replacing by $\delta$ a smaller positive constant if necessary,
we have that for all $v \in D_{\delta \sqrt{n}}$,
\[
e^{n\beta\(\frac{v}{\sqrt{n}}\)}\le \(1-\frac{\pi^2}{n}\abr{v, \SS v}\)^n \le e^{-\pi^2 \abr{v, \SS v}}.
\] 
Since $|\f_{\mu_n}(v/\sqrt{n})| \le C_\delta e^{-\pi^2\abr{v, \SS v}}$ for all $v \in D_{\delta \sqrt{n}}$,
for all large enough $n$,
\begin{align*}
&\frac{1}{n^{\frac{m}{2}}}\int_{D_{\delta \sqrt{n}}\setminus D_{n^{1/8}}}\Big|\f_{\mu_n}\(\frac{v}{\sqrt{n}}\)e^{-2\pi i \abr{v, \F(x)}/\sqrt{n}}\Big|\,dv
\le \frac{C_\delta}{n^{\frac{m}{2}}}\int_{D_{\delta \sqrt{n}}\setminus D_{n^{1/8}}}e^{-\pi^2 \abr{v, \SS v}}\,dv\\
&\le \frac{C_\delta'}{n^{\frac{m}{2}}}\int_{c_0 n^{1/8}}^{c_1\delta \sqrt{n}}e^{-c r^2}r^{m-1}\,dr
\le \frac{C_\delta'}{n^{\frac{m}{2}}}\int_{c_0 n^{1/8}}^{c_1\delta \sqrt{n}}r e^{-\alpha r^2}\,dr
=O\(\frac{1}{n^{\frac{m}{2}}}e^{-\alpha c_0^2n^{1/4}}\)\ll \frac{1}{n^{\frac{m+1}{2}}}.
\end{align*}
In the above, $C_\delta, C_\delta', c, c_0, c_1$ and $\alpha$ (where $c>\alpha$) are positive constants; 
we have used the polar coordinate and the positive definiteness of $\SS$ in the second inequality,
and that $e^{-c r^2}r^{m-1} \le e^{-\alpha r^2}r$ for all large $r$ in the third inequality.

Furthermore for all $n \in \Z_{>0}$ and for all $v \in D_{n^{1/8}}$,
since $\|v\|^4/n \ll 1/\sqrt{n}$,
\[
\f_{\mu_n}\(\frac{v}{\sqrt{n}}\)=e^{-2\pi^2\abr{v, \SS v}}(1+R_n(v)) \(\abr{e^{2\pi i\f/\sqrt{n}}\1, f_{u/\sqrt{n}}}_\pi e^{-2\pi i\f(x_0)/\sqrt{n}} f_{u/\sqrt{n}}(x_0)+O(e^{-c_\delta n})\),
\]
where $R_n(v)=O\(\|v\|^4/n\)$.
For the normalized eigenvector $f_{u}$ of $e^{\beta(v)}$ for $u$ near $0$,
since $f_0=\1$, by Lemma \ref{Lem:der} \eqref{Eq:Lem:eigenvector}
the Taylor theorem shows the following: for each $x \in V(G)$,
\[
f_{u/\sqrt{n}}(x)=1+O\(\frac{\|u\|^2}{n}\).
\]
Furthermore, since $\|u\|\le \|v\|$ and $\|\f\|_\infty \le C\|v\|$, 
for $v \in D_{\delta \sqrt{n}}$,
\begin{align*}
\abr{e^{2\pi i \f/\sqrt{n}}\1, f_{u/\sqrt{n}}}_\pi e^{-2\pi i\f(x_0)/\sqrt{n}}f_{u/\sqrt{n}}(x_0)
&=
\left(1+O\left(\frac{\|\f\|_\infty}{\sqrt{n}}\right)\right)\left(1+O\(\frac{\|v\|^2}{n}\)\right)\\
&=1+O\left(\frac{\|v\|}{\sqrt{n}}\right).
\end{align*}
Note that the change of variable $v\mapsto v/\sqrt{n}$ yields $u\mapsto u/\sqrt{n}$ and $\f\mapsto \f/\sqrt{n}$.
For each $k\ge 0$,
\[
\int_{D_{n^{1/8}}}e^{-2\pi^2\abr{v, \SS v}} \frac{\|v\|^k}{\sqrt{n}}\,dv\le \frac{1}{\sqrt{n}}\int_{\R^m}\|v\|^k e^{-2\pi^2\abr{v, \SS v}}\,dv=O\(\frac{1}{\sqrt{n}}\).
\]
Summarizing the above estimates with $e^{-c_\delta n}\ll n^{-\frac{m+1}{2}}$ yields for $x\in \Gamma$ and for $n \in \Z_{>0}$,
\begin{align}\label{Eq:mu_n_decompose}
\mu_n(x)&=\frac{1}{n^{\frac{m}{2}}}\int_{D_{n^{1/8}}}e^{-2\pi^2\abr{v, \SS v}} e^{-2\pi i \abr{v, \F(x)}/\sqrt{n}}\,dv+O\(\frac{1}{n^{\frac{m+1}{2}}}\).
\end{align}

A direct computation on the Fourier transform yields for all $n \in \Z_{>0}$ and for all $x \in \Gamma$,
\[
\xi_{n\SS}(x)=\frac{1}{n^{\frac{m}{2}}}\int_{\R^m}e^{-2\pi i \abr{v, \F(x)}/\sqrt{n}}e^{-2\pi^2 \abr{v, \SS v}}\,dv.
\]
Abusing notations, we have a constant $\alpha>0$ such that for all $n \in \Z_{>0}$ and uniformly in $x \in \Gamma$,
\begin{equation}\label{Eq:xi_truncation}
\xi_{n\SS}(x)=\frac{1}{n^{\frac{m}{2}}}\int_{D_{n^{1/8}}}e^{-2\pi i \abr{v, \F(x)}/\sqrt{n}}e^{-2\pi^2 \abr{v, \SS v}}\,dv+O\(\frac{1}{n^{\frac{m}{2}}}e^{-\alpha n^{1/4}}\).
\end{equation}
Therefore by \eqref{Eq:mu_n_decompose} and \eqref{Eq:xi_truncation},
for all $n \in \Z_{>0}$ and uniformly in $x \in \Gamma$,
\[
\mu_n(x)=\xi_{n\SS}(x)+O\(\frac{1}{n^{\frac{m+1}{2}}}\).
\]
Furthermore the explicit form of $\SS$ is obtained by Lemma \ref{Lem:H}, as claimed.
\qed

\begin{lemma}\label{Lem:martingale}
Let $\mu$ be a probability measure on $\Gamma$ such that $\supp \mu$ is finite and $\Gamma=\abr{\supp \mu}$, and $\{w_n\}_{n\in \Z_+}$ be a $\mu$-random walk with $w_0=\id$.
There exists a constant $C>0$ such that for all $n \in \Z_{>0}$ and for all real $r >0$,
\[
\Pr\big(|w_n|_S \ge r\big)\le C\exp\Big(-\frac{r^2}{C n}\Big).
\]
\end{lemma}

\proof
The proof follows from the Gaussian estimates established in a more general setting (\cite[the proof of Theorem 9.1]{HebischSaloffCoste93} and \cite[Chapter 14]{WoessBook}).
We provide an alternative proof adapted to this setting for the sake of completeness.

Note that if $\{\F(w_n)\}_{n \in \Z_+}$ is a martingale with respect to the filtration associated with $\{w_n\}_{n \in \Z_+}$, then the proof follows from a concentration inequality.
In general, $\{\F(w_n)\}_{n \in \Z_+}$ is not a martingale.
What we do below is to replace $\F$ by another $\Lambda$-equivariant map which makes the images of $w_n$ form a martingale.

First we claim that there exists a $\Lambda$-{\bf equivariant harmonic map} $\F_H: \Gamma \to \R^m$.
This is a map satisfying the following: $\F_H(gx)=\F_H(x)+g$ for all $x \in \Gamma$ and for all $g \in \Lambda$ 
under the identification between $\Lambda$ and a lattice in $\R^m$,
and
\[
\sum_{s \in \supp \mu}\(\F_H(xs)-\F_H(x)\)\mu(s)=0 \quad \text{for each $x \in \Gamma$}.
\]
This map is obtained from a $\Lambda$-equivariant lift of a Dirichlet energy minimizing map from $G=(V(G), E(G))$ with weights on edges $c(e)$ for $e \in E(G)$ into the flat torus $\R^m/\Lambda$ equipped with metric as a quotient of the standard Euclidean space.
The existence of such a map is shown by a simple variational calculus \cite[Theorem 2.3]{KotaniSunadaStandard}  (see also \cite[Chapter 7]{SunadaTopological}).
(In general, $\F_H$ is not necessarily injective, but this does not affect the following discussion.)
For a $\Lambda$-equivariant harmonic map $\F_H$, we have a martingale $\{\F_H(w_n)\}_{n \in \Z_+}$ with respect to the natural filtration.

Next note that the map $\F_H$ yields a quasi-isometry between $\Cay(\Gamma, S)$ and $\R^m$.
In particular, there exist constants $c_0, c_1>0$ such that 
\begin{equation}\label{Eq:QI}
\|\F_H(x)\|_\infty \ge c_0|x|_S-c_1 \quad \text{for all $x \in \Gamma$}.
\end{equation}
In the above, $\|\cdot\|_\infty$ denotes the $\ell_\infty$-norm in $\R^m$.

Finally, each component of $\F_H(w_n)$ in the coordinate of $\R^m$ is a martingale with a uniformly bounded difference $B$ for some $B>0$.
Hence a union bound and the Azuma-Hoeffding inequality show that by \eqref{Eq:QI},
for all $r \in \Z_+$ and for all $n \in \Z_{>0}$,
\[
\Pr\(|w_n|_S \ge (r+c_1)/c_0\) \le \Pr\(\|\F_H(w_n)\|_\infty \ge r\) \le 2m\exp\Big(-\frac{r^2}{2B^2 n}\Big).
\]
Therefore taking a large enough constant $C>0$ concludes the inequality as claimed.
\qed

\medskip

For the $\xi_{\SS}$ in Theorem \ref{Thm:LCLT},
let us define a {\bf discrete normal distribution} $\Nc^\F_\SS$ on $\Gamma$ by
\[
\Nc^\F_\SS(x):=\frac{1}{Z}\xi_\SS(x) \quad \text{where $Z:=\sum_{x \in \Gamma}\xi_\SS(x)$ for $x \in \Gamma$}.
\]

\begin{theorem}\label{Thm:LCLT_ell1}
In the same setting and assumption as in Theorem \ref{Thm:LCLT},
there exists a constant $C>0$ such that for all integers $n>1$,
\[
\|\mu_n-\Nc^\F_{n\SS}\|_\TV \le \frac{C(\log n)^{\frac{m}{2}}}{\sqrt{n}}.
\]
\end{theorem}

\proof
The local central limit theorem (Theorem \ref{Thm:LCLT}) implies that there exists a constant $C>0$ such that for all $n \in \Z_{>0}$,
\begin{equation}\label{Eq:LCLT}
\sup_{x \in \Gamma}|\mu_n(x)-\xi_{n\SS}(x)|\le \frac{C}{n^{\frac{m+1}{2}}}.
\end{equation}
For a $\mu$-random walk $\{w_n\}_{n \in \Z_+}$ with $w_0=\id$ on $\Gamma$,
Lemma \ref{Lem:martingale}, 
there exists a constant $C>0$ such that for all $n \in \Z_{>0}$ and for all $r>0$,
\begin{equation}\label{Eq:error_martingale}
\Pr\big(|w_n|_S \ge r\big) \le C\exp\Big(-\frac{r^2}{C n}\Big).
\end{equation}
A direct computation on $\xi_{n\SS}$ yields for a (possibly different) constant $C>0$,
for all $n>0$ and for all $r>C\sqrt{n}$,
\begin{equation}\label{Eq:error_Gaussian}
\sum_{|x|_S \ge r}\xi_{n\SS}(x) \le C\exp\Big(-\frac{r^2}{C n}\Big).
\end{equation}
Indeed, this follows from an approximation by a Gaussian $f_{\SS}$ on $\R^m$ for which $\xi_{\SS}=f_\SS\circ \F$, and that $\F$ yields a quasi-isometry between $\Cay(\Gamma, S)$ and $\R^m$.
Note that $\R^m=\bigcup_{x \in \Gamma}x\wbar C_0$ for the closure $\wbar C_0$ of a relatively compact fundamental domain $C_0$ of $\Gamma$.
For $v_i \in \R^m$, $i=1, 2$,
the following holds:
\begin{align*}
|\abr{v_1, \SS^{-1}v_1}-\abr{v_2, \SS^{-1}v_2}|
&= \left|\int_0^1 \frac{d}{dt}\abr{v_1+t(v_2-v_1), \SS^{-1}(v_1+t(v_2-v_1))}\,dt\right| \\
&\le 2\int_0^1 |\abr{v_2-v_1, \SS^{-1}(v_1+t(v_2-v_1))}|\,dt \\
&\le 2\|\SS^{-1}\|\|v_1-v_2\|\max\{\|v_1\|, \|v_2\|\}.
\end{align*}
Letting $\diam C_0$ denote the diameter of $\wbar C_0$,
we have that by the above inequality, if $v_i \in x\wbar C_0$ and $\|v_i\| \ge \diam C_0$, $i=1, 2$, 
then $\|v_1-v_2\| \le \|v_2\|$ and
\[
|\abr{v_1, \SS^{-1}v_1}-\abr{v_2, \SS^{-1}v_2}|\le 2\|\SS^{-1}\|\|v_1-v_2\|(\|v_2\|+\|v_1-v_2\|)
\le 4\|\SS^{-1}\|\diam C_0 \|v_2\|.
\]
Thus, for all $v_i \in x \wbar C_0$, $i=1, 2$,
\[
f_{n\SS}(v_1) \le f_{n\SS}(v_2)e^{\frac{c}{n}\|v_2\|}, \quad \text{where $c:=2\|\SS^{-1}\|\diam C_0$}.
\]
Since $\SS^{-1}$ is positive definite, there exists a constant $\alpha>0$ such that
$\abr{v, \SS^{-1}v} \ge \alpha\|v\|^2$ for $v \in \R^m$.
Therefore noting that $\|\F(x)\| \ge c_0|x|_S-c_1$ for all $x \in \Gamma$,
we estimate
\[
\sum_{|x|_S \ge r} \xi_{n\SS}(x) \le \frac{1}{{\rm vol}(C_0)}\int_{\|v\| \ge c_0 r-c_1}\frac{1}{(2\pi n)^{\frac{m}{2}}\sqrt{\det \SS}}e^{-\frac{\alpha}{2n}\|v\|^2+\frac{c}{n}\|v\|}\,dv,
\]
for $r \ge (\diam C_0+c_1)/c_0$, where ${\rm vol}(C_0)$ is the volume of $C_0$.
By the change of variables $v \mapsto \sqrt{n}v$, the right hand side equals the following:
\begin{align*}
&\frac{1}{{\rm vol}(C_0)}\int_{\|v\| \ge (c_0 r-c_1)/\sqrt{n}}\frac{1}{(2\pi)^{\frac{m}{2}}\sqrt{\det \SS}}e^{-\frac{\alpha}{2}\|v\|^2+\frac{c}{\sqrt{n}}\|v\|}\,dv\\
&\qquad \qquad \qquad \qquad \qquad \qquad \qquad \qquad \qquad
=c_{m, \SS}\int_{(c_0 r-c_1)/\sqrt{n}}^\infty e^{-\frac{\alpha}{2}s^2+\frac{c}{\sqrt{n}}s}s^{m-1}\,ds,
\end{align*}
where $c_{m, \SS}$ is a constant depending only on $\SS$ and $m$ in the polar coordinate.
Thus,
\begin{align*}
\int_{(c_0 r-c_1)/\sqrt{n}}^\infty e^{-\frac{\alpha}{2}\(s-\frac{c}{\alpha\sqrt{n}}\)^2+\frac{c^2}{2\alpha n}}s^{m-1}\,ds
=e^{\frac{c^2}{2\alpha n}}\int_{R}^\infty e^{-\frac{\alpha}{2}s^2}\(s+\frac{c}{\alpha \sqrt{n}}\)^{m-1}\,ds,
\end{align*}
where $R:=(c_0 r-c_1)/\sqrt{n}-c/(\alpha \sqrt{n})$ by change of variables.
There exists a constant $C>0$ such that 
$\(s+c/(\alpha \sqrt{n})\)^{m-1}\le s e^{\frac{\alpha}{4}s^2}$ for all $s>C$.
Hence there exists a constant $C>0$ such that for all $n>0$ and $r> C\sqrt{n}$,
the last term is at most
\[
e^{\frac{c^2}{2\alpha n}}\int_{R}^\infty se^{-\frac{\alpha}{4}s^2}\,ds
=e^{\frac{c^2}{2\alpha n}}\frac{2}{\alpha}e^{-\frac{\alpha}{4}R^2}
\le C e^{-\frac{r^2}{C n}}.
\]
This shows \eqref{Eq:error_Gaussian}.

Combining \eqref{Eq:LCLT}, \eqref{Eq:error_martingale} and \eqref{Eq:error_Gaussian} yields for all $n>0$ and for all $r>C\sqrt{n}$,
\begin{align*}
\|\mu_n-\xi_{n\SS}\|_1	&=\sum_{|x|_S \le r}|\mu_n(x)-\xi_{n\SS}(x)|+\sum_{|x|_S > r}|\mu_n(x)-\xi_{n\SS}(x)|\\
						&\le |B_{S}(r)|\frac{C}{n^{\frac{m+1}{2}}}+2C\exp\Big(-\frac{r^2}{C n}\Big).
\end{align*}
In the above $B_S(r):=\{x \in \Gamma \ : \ |x|_S \le r\}$.
Note that $|B_S(r)|=\Theta(r^m)$, in particular, $|B_S(r)|\le C r^m$ for all $r \in \Z_+$ for a constant $C>0$ independent of $r$ since $\F$ is a $\Lambda$-equivariant injective map from $\Gamma$ into $\R^m$.
Letting $r=A\sqrt{n \log n}$ for a constant $A>0$, we have $r>C\sqrt{n}$ for all $n \in \Z_{>0}$,
and
\[
\|\mu_n-\xi_{n\SS}\|_1 \le \frac{C^2 A^m(\log n)^{\frac{m}{2}}}{\sqrt{n}}+\frac{2C}{n^{A^2/C}}.
\]
Fixing a large enough constant $A$ such that $A^2/C>1/2$ shows that there exists a constant $C_1$ such that for all $n>1$,
\begin{equation}\label{Eq:ell1}
\|\mu_n-\xi_{n\SS}\|_1 \le \frac{C_1(\log n)^{\frac{m}{2}}}{\sqrt{n}}.
\end{equation}
Note that there exists a constant $c>0$ such that for all $n \in \Z_{>0}$,
\begin{equation}\label{Eq:discrete_normal}
\sum_{x \in \Gamma}\xi_{n\SS}(x)=1+\sum_{x \in \Gamma\setminus\{\id\}}\xi_{n\SS}(x)=1+O(e^{-cn}).
\end{equation}
Indeed, since $\Z^m$ is a finite index subgroup of $\Gamma$ and $\F$ is a $\Z^m$-equivariant embedding with a discrete image in $\R^m$, 
the Poisson summation formula
(cf.\ \eqref{Eq:Poisson} in Appendix \ref{Sec:appendix}) on finitely many orbits of $\Z^m$ shows \eqref{Eq:discrete_normal}.
This implies that for all $n \in \Z_{>0}$,
\[
\|\xi_{n\SS}-\Nc^\F_{n\SS}\|_1 =\sum_{x \in \Gamma}\Big|\xi_{n\SS}(x)-\frac{1}{\sum_{z \in \Gamma}\xi_{n\SS}(z)}\xi_{n\SS}(x)\Big|=O(e^{-c n}).
\]
Therefore this together with \eqref{Eq:ell1}, adjusting a constant factor $C$ yields the claim.
\qed

\begin{remark}\label{Rem:sharp}
The inequality \eqref{Eq:ell1} suffices for our purpose in Theorem \ref{Thm:NSonWeyl} below.
If $\Phi$ is $\mu$-harmonic, then the estimate in Theorem \ref{Thm:LCLT_ell1} becomes $O((\log n)^{\frac{m}{2}}/n)$.
In this case, the estimate is sharp up the factor $O((\log n)^{m/2})$. 
Indeed, the local central limit theorem provides an example
satisfying that
$\|\mu_n-\Nc^\F_{n\SS}\|_\TV=\Omega(1/n)$,
e.g., a lazy simple random walk on $\Z^m$.
\end{remark}

\section{Applications to noise sensitivity problem}\label{Sec:applications}

\subsection{Noise sensitivity on affine Weyl groups}\label{Sec:NSonWeyl}

Let $(\Gamma, S)$ be an affine Weyl group where $S$ is a canonical set of generators consisting of involutions.
We have $\Gamma=\Lambda \rtimes W$ as in Section \ref{Sec:affineWeyl}.
For a probability measure $\mu$ on $\Gamma$ and for all $\rho \in [0, 1]$,
we recall that $\pi^\rho=\rho (\mu\times \mu)+(1-\rho)\mu_\diag$,
where $\mu \times \mu$ denote the product measure and $\mu_\diag((x, y))=\mu(x)$ if $x=y$ and $0$ otherwise.

For each $\rho \in [0, 1]$, the measure $\pi^\rho$ is defined on $\Gamma \times \Gamma=(\Lambda\times \Lambda)\rtimes (W\times W)$, for which $S_\ast:=(S\cup \{\id\})^2$ is a generating set of order at most $2$.
On the one hand, if $\supp \mu=S\cup \{\id\}$,
then $\supp \pi^\rho=S_\ast$ for all $\rho \in (0, 1]$.
On the other hand, however, if $\rho=0$ or $\supp \mu=S$ (in which case, the support of $\mu$ does not contain $\id$), then $\supp \pi^\rho$ {\em never} generate $\Gamma \times \Gamma$.
Indeed, every affine Weyl group $\Gamma$ admits a surjective homomorphism onto $\{\pm 1\}$ through the determinant of the isometry part in the natural affine representation.
The product group $\Gamma \times \Gamma$ admits a surjective homomorphism onto $\{\pm 1\}^2$.
If $\supp \mu$ does not contain $\id$, then $\supp \pi^\rho$ for $\rho \in (0, 1]$ only generates a proper subgroup of $\Gamma \times \Gamma$ (of index $2$).
Furthermore if $\rho=0$, then $\supp \pi^\rho$ generates the diagonal subgroup isomorphic to $\Gamma$ in $\Gamma \times \Gamma$.

\begin{theorem}\label{Thm:NSonWeyl}
Let $(\Gamma, S)$ be an affine Weyl group, and $\mu$ be a probability measure on $\Gamma$ such that
the support of $\mu$ equals $S\cup\{\id\}$.
For all $\rho \in (0, 1]$, the $\pi^\rho$-random walk on $\Gamma\times \Gamma$ starting from the identity satisfies the following:
There exist a constant $C>0$ and an integer $m>0$ such that for all integers $n>1$,
\[
\|\pi^\rho_n-\mu_n\times \mu_n\|_\TV \le \frac{C(\log n)^m}{\sqrt{n}}.
\]
In particular, for all $\rho \in (0, 1]$,
\[
\lim_{n \to \infty}\|\pi^\rho_n-\mu_n\times \mu_n\|_\TV=0,
\]
i.e., the $\mu$-random walk on $\Gamma$ is noise sensitive in total variation.
\end{theorem}

\proof
Let us apply to $\Gamma \times \Gamma$ and $\pi^\rho$ the discussion we have made so far.
Fix a point $o$ in the interior of a chamber for $\Gamma$ in the associated Euclidean space $\R^m$.
We assume that $o$ is the origin after taking a conjugate by a translation for the action of $\Gamma$ if necessary.
Let $\F:=\F^{(1)}\times \F^{(2)}: \Gamma \times \Gamma \to \R^m \times \R^m$, $(x_1, x_2)\mapsto (x_1.o, x_2.o)$,
where $\F^{(i)}(x)=x.o$ for $x \in \Gamma$ and for $i=1, 2$.
For all $\rho \in (0, 1]$, the support of $\pi^\rho$ is finite, $\Gamma\times \Gamma=\abr{\supp \pi^\rho}$, 
and $\pi^\rho$ is symmetric since every element in $\supp \pi^\rho=S_\ast=(S\cup\{\id\})^2$ has the order at most $2$.
Further let $G:=(\Lambda \times \Lambda)\backslash \Cay(\Gamma\times \Gamma, S_\ast)$ equipped with the conductance $c(\eb)$ for $\eb \in E(G)$ induced from $\pi^\rho$.
We consider the corresponding pointed finite network $(G, c, x_0)$ where $x_0$ is the identity element in $W\times W=(\Lambda\times \Lambda) \backslash (\Gamma \times \Gamma)$.
For all $\rho \in (0, 1]$, the corresponding Markov chain on $G$ is irreducible and satisfies that $p(x, x)>0$ for every $x \in V(G)$ since $\pi^\rho(\id)>0$.

Let $\SS^\rho$ be the matrix for $\pi^\rho$ in the local central limit theorem (Theorem \ref{Thm:LCLT}),
let us show that $\SS^\rho=\SS^1$ for all $\rho \in (0, 1]$.
Namely, in the block diagonal form along the decomposition $\R^m\times \R^m$, 
we prove the following: 
For all $\rho \in (0, 1]$,
\[
\SS^\rho=
\begin{pmatrix}
\SS^\mu & 0\\
0 & \SS^\mu
\end{pmatrix},
\]
where $\SS^\mu$ is the matrix $\SS^\mu$ for $\mu$ in the local central limit theorem.
Theorem \ref{Thm:LCLT} shows that $\SS^\rho$ is obtained by
\begin{equation}\label{Eq:Thm:covariance}
\abr{\vb_1, \SS^\rho \vb_2}=\sum_{\eb \in E(G)}\ub_1(\eb)\ub_2(\eb)c(\eb),
\end{equation}
for $\vb_1, \vb_2 \in \R^m\times\R^m$ and $\ub_i$ is the harmonic part of the $1$-form $\wh \vb_i=(\abr{\vb_i, \F_\eb})_{\eb \in E(G)}$ for $i=1, 2$.

First let us show that $\abr{\vb_1, \SS^\rho \vb_2}=0$ for all $\vb_1=(v_1, 0), \vb_2=(0, v_2) \in \R^m\times \R^m$.
For each edge $\eb=((x_1, x_2), (x_1.s_1, x_2.s_2)) \in E(G)$, 
we write $e_i=(x_i, x_i.s_i)$ for $x_i \in W$ and $s_i \in S\cup\{\id\}$ for $i=1, 2$.
For such $\vb_1, \vb_2$, we have that
\begin{equation}\label{Eq:v}
\wh \vb_i(\eb)=\abr{\vb_i, \F_\eb}=\abr{v_i, \F^{(i)}_{e_i}}=\wh v_i(e_i) \quad \text{for $i=1, 2$}.
\end{equation}
Let $\ub_i$ be the harmonic part of $\wh \vb_i$, and $u_i$ be the harmonic part of $\wh v_i$ for each $i=1, 2$.
It holds that
\begin{equation}\label{Eq:harmonic_part}
\ub_i(\eb)=u_i(e_i) \quad \text{for $i=1, 2$}.
\end{equation}
Indeed, let $\wt u_i(\eb):=u_i(e_i)$ for $i=1, 2$.
These define harmonic $1$-forms on $G$: 
This follows since $S_\ast=(S\cup\{\id\})^2$ and each marginal of $\pi^\rho$ is $\mu$.
Furthermore by the definition of $u_i$,
it holds that for some $f_i:W \to \R$, 
\begin{equation}\label{Eq:f}
\wh v_i(x, x.s)=u_i(x, x.s)+df_i(x, x.s) \quad \text{for $x \in W$ and $s \in S$}.
\end{equation}
Hence by \eqref{Eq:v} and \eqref{Eq:f}, letting $\wt f_i: V(G) \to \R$ by $\wt f_i((x_1, x_2)):=f_i(x_i)$ for $(x_1, x_2)\in V(G)=W\times W$ for $i=1, 2$ yields
\[
\wh \vb_i(\eb)=\wt u_i(\eb)+d\wt f_i(\eb) \quad \text{for $\eb \in E(G)$}.
\]
The uniqueness of harmonic part concludes \eqref{Eq:harmonic_part}.

The right hand side in \eqref{Eq:Thm:covariance} is computed as in the following: 
By \eqref{Eq:harmonic_part}, 
it holds that
\begin{align*}
&\sum_{\eb \in E(G)}\ub_1(\eb)\ub_2(\eb)c(\eb)
=\sum_{\eb \in E(G)}u_1(e_1)u_2(e_2)c(\eb)\\
&=\sum_{(x_1, x_2)\in W\times W}\sum_{(s_1, s_2) \in S_\ast}u_1(x_1, x_1.s_1)u_2(x_2, x_2.s_2)\pi((x_1, x_2))\pi^\rho((s_1, s_2)).
\end{align*}
Note that the summation over $(s_1, s_2) \in S_\ast$ is restricted to $S \times S$
since $u_i(x, x)=0$ for $x \in W$ and for each $i=1, 2$.
Furthermore $\pi((x, y))=1/|W|^2$ for all $(x, y) \in W\times W$.
Hence the last term times the factor $|W|^2$ leads the following: 
\begin{align*}
&\sum_{x_1 \in W, s_1 \in S}u_1(x_1, x_1.s_1)\sum_{x_2 \in W, s_2 \in S}u_2(x_2, x_2.s_2)\pi^\rho((s_1, s_2))\\
&=\sum_{x_1 \in W, s_1 \in S}u_1(x_1, x_1.s_1)\frac{1}{2}\sum_{x_2 \in W, s_2 \in S}(u_2(x_2, x_2.s_2)+u_2(x_2.s_2, x_2))\pi^\rho((s_1, s_2))=0.
\end{align*}
In the above we have used that $s_2=s_2^{-1}$ and $\pi^\rho((s_1, s_2))=\pi^\rho((s_1, s_2^{-1}))$ for each $(s_1, s_2) \in S\times S$ in the first equality, and that $u_2(x_2, x_2.s_2)=-u_2(x_2.s_2, x_2)$ for all $x_2 \in W$ and all $s_2 \in S$ in the last equality.
Therefore for all $\vb_1=(v_1, 0), \vb_2=(0, v_2) \in \R^m \times \R^m$,
it holds that
$\abr{\vb_1, \SS^\rho \vb_2}=0$.

Next since each marginal of $\pi^\rho$ is $\mu$,
for $\vb_1=(v_1, 0), \vb_2=(v_2, 0) \in \R^m \times \R^m$,
\begin{align*}
\abr{\vb_1, \SS^\rho \vb_2}
=\abr{v_1, \SS^\mu v_2}.
\end{align*}
The same equality holds for $\vb_1=(0, v_1), \vb_2=(0, v_2) \in \R^m \times \R^m$.
Summarizing all the above discussion, we obtain $\SS^\rho=\SS^1$ for all $\rho \in (0, 1]$.

Finally, Theorem \ref{Thm:LCLT_ell1} shows that for all $\rho \in (0, 1]$,
there exists a constant $C>0$ such that for all $n>1$, by the triangle inequality, 
\[
\|\pi^\rho_n-\mu_n\times \mu_n\|_\TV \le \|\pi^\rho_n-\Nc^\F_{n\SS^\rho}\|_\TV+\|\mu_n\times \mu_n-\Nc^\F_{n\SS^1}\|_\TV\le \frac{2C(\log n)^m}{\sqrt{n}}.
\]
This concludes the first claim. 
The second claim follows from the first claim.
\qed

\subsection{Examples}\label{Sec:example}

Let us provide explicit examples of random walks on the affine Weyl groups of type $\wt A_1 \times \wt A_1$, $\wt A_2$ and $\wt C_2$.
Figures \ref{Fig:A1A1}, \ref{Fig:A2} and \ref{Fig:C2} respectively describe
the Cayley graphs of the groups with the corresponding sets of generators (the solid lines with dots). 
The associated action of each group on $\R^2$ consists of reflections with respect to lines (indicated as dotted lines) with a fundamental domain (colored in dark gray). The lattice has a larger fundamental domain (colored in light gray).

\subsubsection{Type $\wt A_1\times \wt A_1$}\label{Sec:exampleA1A1}

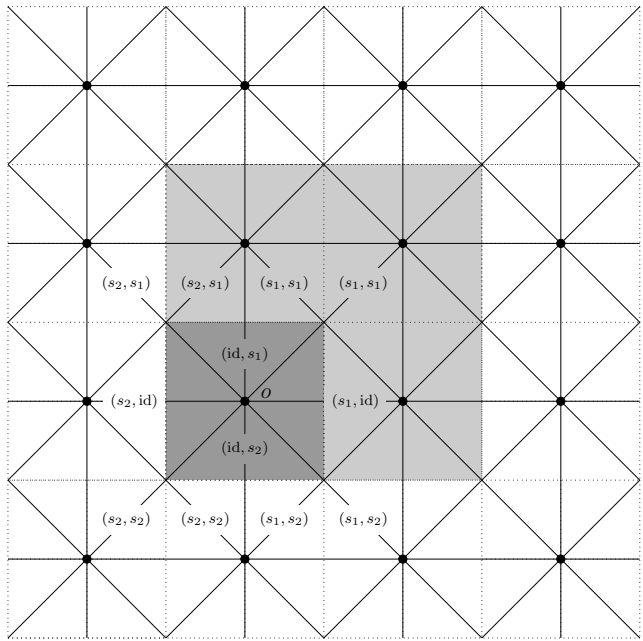
\begin{figure}[h]
\centering
\scalebox{0.7}[0.7]{
\begin{tikzpicture}[
dot/.style={draw, coordinate},
ver/.style={draw, circle, scale=0.4, fill=black},
]



\node[dot] at (-4.5,7.5) (-4.5,7.5) {};
\node[dot] at (-1.5,7.5) (-1.5,7.5) {};
\node[dot] at (1.5,7.5) (1.5,7.5) {};
\node[dot] at (4.5,7.5) (4.5,7.5) {};
\node[dot] at (7.5,7.5) (7.5,7.5) {};

\node[dot] at (-4.5,4.5) (-4.5,4.5) {};
\node[dot] at (-1.5,4.5) (-1.5,4.5) {};
\node[dot] at (1.5,4.5) (1.5,4.5) {};
\node[dot] at (4.5,4.5) (4.5,4.5) {};
\node[dot] at (7.5,4.5) (7.5,4.5) {};

\node[dot] at (-4.5,1.5) (-4.5,1.5) {};
\node[dot] at (-1.5,1.5) (-1.5,1.5) {};
\node[dot] at (1.5,1.5) (1.5,1.5) {};
\node[dot] at (4.5,1.5) (4.5,1.5) {};
\node[dot] at (7.5,1.5) (7.5,1.5) {};

\node[dot] at (-4.5,-1.5) (-4.5,-1.5) {};
\node[dot] at (-1.5,-1.5) (-1.5,-1.5) {};
\node[dot] at (1.5,-1.5) (1.5,-1.5) {};
\node[dot] at (4.5,-1.5) (4.5,-1.5) {};
\node[dot] at (7.5,-1.5) (7.5,-1.5) {};

\draw[dotted, fill=gray!40] (-1.5,-1.5)--(4.5,-1.5)--(4.5,4.5)--(-1.5,4.5)--cycle;
\draw[dotted, fill=gray!80] (-1.5,-1.5)--(1.5,-1.5)--(1.5,1.5)--(-1.5,1.5)--cycle;

\node[dot] at (-4.5,-4.5) (-4.5,-4.5) {};
\node[dot] at (-1.5,-4.5) (-1.5,-4.5) {};
\node[dot] at (1.5,-4.5) (1.5,-4.5) {};
\node[dot] at (4.5,-4.5) (4.5,-4.5) {};
\node[dot] at (7.5,-4.5) (7.5,-4.5) {};

\node[dot] at (-4.5,-7.5) (-4.5,-7.5) {};
\node[dot] at (-1.5,-7.5) (-1.5,-7.5) {};
\node[dot] at (1.5,-7.5) (1.5,-7.5) {};
\node[dot] at (4.5,-7.5) (4.5,-7.5) {};
\node[dot] at (7.5,-7.5) (7.5,-7.5) {};


\draw[dotted] (4.5,3)--(7.5,3);
\draw[dotted] (4.5,0)--(7.5,0);
\draw[dotted] (4.5,-3)--(7.5,-3);
\draw[dotted] (-4.5,6)--(-1.5,6);
\draw[dotted] (-4.5,3)--(-1.5,3);
\draw[dotted] (-4.5,0)--(-1.5,0);
\draw[dotted] (-4.5,-3)--(-1.5,-3);

\draw[dotted] (-4.5,7.5)--(-1.5,7.5);
\draw[dotted] (1.5,7.5)--(-1.5,7.5);
\draw[dotted] (1.5,7.5)--(4.5,7.5);
\draw[dotted] (7.5,7.5)--(4.5,7.5);

\draw[dotted] (-4.5,4.5)--(-1.5,4.5);
\draw[dotted] (1.5,4.5)--(-1.5,4.5);
\draw[dotted] (1.5,4.5)--(4.5,4.5);
\draw[dotted] (7.5,4.5)--(4.5,4.5);

\draw[dotted] (-1.5,1.5)--(1.5,1.5);
\draw[dotted] (4.5,1.5)--(1.5,1.5);
\draw[dotted] (-1.5,4.5)--(1.5,4.5);
\draw[dotted] (4.5,1.5)--(7.5,1.5);
\draw[dotted] (-4.5,1.5)--(-1.5,1.5);
\draw[dotted] (4.5,4.5)--(7.5,4.5);

\draw[dotted] (-1.5,-1.5)--(-4.5,-1.5);
\draw[dotted] (4.5,-1.5)--(7.5,-1.5);
\draw[dotted] (-4.5,-1.5)--(-1.5,-1.5);
\draw[dotted] (1.5,-1.5)--(4.5,-1.5);

\draw[dotted] (-4.5,-4.5)--(-1.5,-4.5);
\draw[dotted] (1.5,-4.5)--(-1.5,-4.5);
\draw[dotted] (-4.5,-4.5)--(-1.5,-4.5);
\draw[dotted] (7.5,-4.5)--(4.5,-4.5);
\draw[dotted] (1.5,-4.5)--(4.5,-4.5);


\draw[dotted] (-4.5,7.5)--(-4.5,4.5);
\draw[dotted] (1.5,7.5)--(1.5,4.5);
\draw[dotted] (-1.5,7.5)--(-1.5,4.5);
\draw[dotted] (4.5,7.5)--(4.5,4.5);
\draw[dotted] (7.5,7.5)--(7.5,4.5);

\draw[dotted] (-4.5,1.5)--(-4.5,4.5);
\draw[dotted] (1.5,1.5)--(1.5,4.5);
\draw[dotted] (-1.5,1.5)--(-1.5,4.5);
\draw[dotted] (4.5,1.5)--(4.5,4.5);
\draw[dotted] (7.5,1.5)--(7.5,4.5);

\draw[dotted] (-4.5,-1.5)--(-4.5,-4.5);
\draw[dotted] (1.5,-1.5)--(1.5,-4.5);
\draw[dotted] (-1.5,-1.5)--(-1.5,-4.5);
\draw[dotted] (4.5,-1.5)--(4.5,-4.5);
\draw[dotted] (7.5,-1.5)--(7.5,-4.5);

\draw[dotted] (-3,-1.5)--(-3,-4.5);
\draw[dotted] (0,-1.5)--(0,-4.5);
\draw[dotted] (3,-1.5)--(3,-4.5);
\draw[dotted] (6,-1.5)--(6,-4.5);

\draw[dotted] (1.5,1.5)--(1.5,-1.5);
\draw[dotted] (4.5,1.5)--(4.5,-1.5);
\draw[dotted] (7.5,1.5)--(7.5,-1.5);
\draw[dotted] (-1.5,1.5)--(-1.5,-1.5);
\draw[dotted] (-4.5,1.5)--(-4.5,-1.5);

\draw[dotted] (0,4.5)--(0,7.5);
\draw[dotted] (3,4.5)--(3,7.5);
\draw[dotted] (1.5,4.5)--(1.5,7.5);
\draw[dotted] (4.5,4.5)--(4.5,7.5);
\draw[dotted] (-1.5,4.5)--(-1.5,7.5);


\node[ver] at (-3,6) (-3,6) {};
\node[ver] at (0,6) (0,6) {};
\node[ver] at (3,6) (3,6) {};
\node[ver] at (6,6) (6,6) {};

\node[ver] at (-3,3) (-3,3) {};
\node[ver] at (0,3) (0,3) {};
\node[ver] at (3,3) (3,3) {};
\node[ver] at (6,3) (6,3) {};

\node[ver] at (-3,0) (-3,0) {};
\node[ver] at (0,0) (0,0) {};
\node[] at (0.4,0.15) (0.4,0.15) {$o$};
\node[ver] at (3,0) (3,0) {};
\node[ver] at (6,0) (6,0) {};

\node[ver] at (-3,-3) (-3,-3) {};
\node[ver] at (0,-3) (0,-3) {};
\node[ver] at (3,-3) (3,-3) {};
\node[ver] at (6,-3) (6,-3) {};


\draw[] (-3,6) -- (-4.5,6);
\draw[] (-3,6) -- (0,6);
\draw[] (0,6) -- (3,6);
\draw[] (3,6) -- (6,6);
\draw[] (7.5,6) -- (6,6);

\draw[] (-3,3) -- (-4.5,3);
\draw[] (-3,3) -- (0,3);
\draw[] (0,3) -- (3,3);
\draw[] (3,3) -- (6,3);
\draw[] (7.5,3) -- (6,3);

\draw[] (-3,0) -- (-4.5,0);
\draw[] (-3,0) -- (0,0) node[pos=0.3, fill=white] {{\tiny$(s_2, \id)$}};
\draw[] (0,0) -- (3,0) node[pos=0.7, fill=gray!40] {{\tiny$(s_1, \id)$}};
\draw[] (3,0) -- (6,0);
\draw[] (7.5,0) -- (6,0);

\draw[] (-3,-3) -- (-4.5,-3);
\draw[] (-3,-3) -- (0,-3);
\draw[] (0,-3) -- (3,-3);
\draw[] (3,-3) -- (6,-3);
\draw[] (7.5,-3) -- (6,-3);


\draw[] (-3,6) -- (-3,7.5);
\draw[] (0,6) -- (0,7.5);
\draw[] (3,6) -- (3,7.5);
\draw[] (6,6) -- (6,7.5);

\draw[] (-3,6) -- (-3,3);
\draw[] (0,6) -- (0,3);
\draw[] (3,6) -- (3,3);
\draw[] (6,6) -- (6,3);

\draw[] (-3,0) -- (-3,3);
\draw[] (0,0) -- (0,3) node[pos=0.3, fill=gray!80] {{\tiny$(\id, s_1)$}};
\draw[] (3,0) -- (3,3);
\draw[] (6,0) -- (6,3);

\draw[] (-3,0) -- (-3,-3);
\draw[] (0,0) -- (0,-3) node[pos=0.3, fill=gray!80] {{\tiny$(\id, s_2)$}};
\draw[] (3,0) -- (3,-3);
\draw[] (6,0) -- (6,-3);

\draw[] (-3,-4.5) -- (-3,-3);
\draw[] (0,-4.5) -- (0,-3);
\draw[] (3,-4.5) -- (3,-3);
\draw[] (6,-4.5) -- (6,-3);


\draw[] (-4.5,-4.5) -- (-3,-3);
\draw[] (-1.5,-4.5) -- (-3,-3);
\draw[] (1.5,-4.5) -- (0,-3);
\draw[] (-1.5,-4.5) -- (0,-3);
\draw[] (1.5,-4.5) -- (3,-3);
\draw[] (4.5,-4.5) -- (3,-3);
\draw[] (4.5,-4.5) -- (6,-3);
\draw[] (7.5,-4.5) -- (6,-3);

\draw[] (0,0) -- (-3,-3) node[pos=0.75, fill=white] {{\tiny$(s_2, s_2)$}};
\draw[] (0,0) -- (3,-3) node[pos=0.75, fill=white] {{\tiny$(s_1, s_2)$}};
\draw[] (6,0) -- (3,-3);
\draw[] (3,0) -- (6,-3);
\draw[] (3,0) -- (0,-3) node[pos=0.75, fill=white] {{\tiny$(s_1, s_2)$}};
\draw[] (-3,0) -- (0,-3) node[pos=0.75, fill=white] {{\tiny$(s_2, s_2)$}};

\draw[] (0,0) -- (-3,3) node[pos=0.75, fill=white] {{\tiny$(s_2, s_1)$}};
\draw[] (0,0) -- (3,3) node[pos=0.75, fill=gray!40] {{\tiny$(s_1, s_1)$}};
\draw[] (6,0) -- (3,3);
\draw[] (3,0) -- (6,3);
\draw[] (3,0) -- (0,3) node[pos=0.75, fill=gray!40] {{\tiny$(s_1, s_1)$}};
\draw[] (-3,0) -- (0,3) node[pos=0.75, fill=gray!40] {{\tiny$(s_2, s_1)$}};

\draw[] (-4.5,4.5) -- (-3,3);
\draw[] (0,6) -- (-3,3);
\draw[] (0,6) -- (3,3);
\draw[] (6,6) -- (3,3);
\draw[] (6,3) -- (7.5,4.5);
\draw[] (6,3) -- (7.5,1.5);
\draw[] (6,0) -- (7.5,1.5);
\draw[] (6,0) -- (7.5,-1.5);
\draw[] (6,-3) -- (7.5,-1.5);
\draw[] (6,6) -- (7.5,4.5);
\draw[] (3,6) -- (6,3);
\draw[] (3,6) -- (0,3);
\draw[] (-3,6) -- (0,3);
\draw[] (-3,6) -- (-4.5,4.5);
\draw[] (-3,3) -- (-4.5,1.5);
\draw[] (-3,0) -- (-4.5,1.5);
\draw[] (-3,0) -- (-4.5,-1.5);
\draw[] (-3,-3) -- (-4.5,-1.5);

\draw[] (-4.5,7.5) -- (-3,6);
\draw[] (-1.5,7.5) -- (-3,6);
\draw[] (1.5,7.5) -- (0,6);
\draw[] (-1.5,7.5) -- (0,6);
\draw[] (1.5,7.5) -- (3,6);
\draw[] (4.5,7.5) -- (3,6);
\draw[] (4.5,7.5) -- (6,6);
\draw[] (7.5,7.5) -- (6,6);

\end{tikzpicture}
}
\caption{The Cayley graph of $D_\infty \times D_\infty$ with the set of generators $S_\ast$ (where loops corresponding to $\id$ are omitted) and 
the group action on $\R^2$.}
\label{Fig:A1A1}
\end{figure}

Let us consider the infinite dihedral group:
\[
D_\infty=\abr{s_1, s_2 \mid s_1^2=s_2^2=\id}.
\]
Let $S:=\{s_1, s_2\}$. 
The pair $(D_\infty, S)$ is the affine Weyl group of type $\wt A_1$, and
the product group $D_\infty \times D_\infty$ with the standard set of generators $S\times\{\id\}\cup\{\id\}\times S$ is the affine Weyl group of type $\wt A_1 \times \wt A_1$.
The group $D_\infty \times D_\infty$ is isomorphic to $\Z^2 \rtimes (\Z/2)^2$.
We define $\F=\F^{(1)}\times \F^{(2)}:D_\infty\times D_\infty \to \R^2$ in a way that the origin is the barycenter of a fundamental chamber (which is a square of side length $1/2$).
The lattice is identified with the standard integer lattice (Figure \ref{Fig:A1A1}).

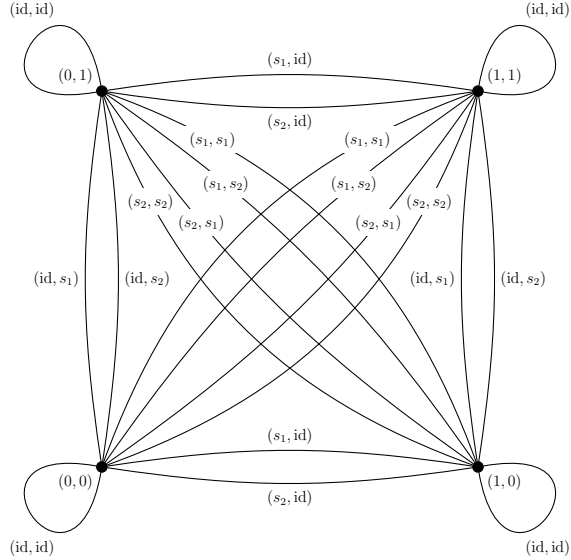
\begin{figure}[h]
\centering
\scalebox{0.5}[0.5]{
\begin{tikzpicture}
\begin{scope}[every node/.style={circle, fill=black, draw, inner sep=0pt, minimum size=8pt}]
\node (AA) at (0, 0) {};
\node (BA) at (10, 0) {};
\node (AB) at (0, 10) {};
\node (BB) at (10, 10) {};
\end{scope}

\begin{scope}[every node/.style={}]
\node at (-0.7, -0.4) {$(0,0)$};
\node at (10.7, -0.4) {$(1,0)$};
\node at (-0.7, 10.4) {$(0,1)$};
\node at (10.7, 10.4) {$(1,1)$};
\end{scope}

\begin{scope}[>={Stealth[black]},
				every node/.style={fill=white},
				every edge/.style={draw=black, thick},
				every loop/.style={}]

	\path [-] (AA) edge[bend left=8.5] node[pos=0.5, above=0.3mm] {$(s_1, {\rm id})$} (BA);
	\path [-] (AA) edge[bend right=8.5] node[pos=0.5, below=0.3mm] {$(s_2, {\rm id})$} (BA);
	
	\path [-] (AB) edge[bend left=8.5] node[pos=0.5, above=0.3mm] {$(s_1, {\rm id})$} (BB);
	\path [-] (AB) edge[bend right=8.5] node[pos=0.5, below=0.3mm] {$(s_2, {\rm id})$} (BB);
	
	\path [-] (AA) edge[bend left=8.5] node[pos=0.5, left=0.1mm] {$({\rm id}, s_1)$} (AB);
	\path [-] (AA) edge[bend right=8.5] node[pos=0.5, right=0.1mm] {$({\rm id}, s_2)$} (AB);
	
	\path [-] (BA) edge[bend left=8.5] node[pos=0.5, left=0.1mm] {$({\rm id}, s_1)$} (BB);
	\path [-] (BA) edge[bend right=8.5] node[pos=0.5, right=0.1mm] {$({\rm id}, s_2)$} (BB);

	\path [-] (AA) edge[bend left=25] node[pos=0.8] {$(s_1, s_1)$} (BB);
	\path [-] (AA) edge[bend right=25] node[pos=0.8] {$(s_2, s_2)$} (BB);
	\path [-] (AA) edge[bend left=10] node[pos=0.73] {$(s_1, s_2)$} (BB);
	\path [-] (AA) edge[bend right=10] node[pos=0.71] {$(s_2, s_1)$} (BB);
	
	\path [-] (AB) edge[bend left=25] node[pos=0.2] {$(s_1, s_1)$} (BA);
	\path [-] (AB) edge[bend right=25] node[pos=0.2] {$(s_2, s_2)$} (BA);
	\path [-] (AB) edge[bend left=10] node[pos=0.27] {$(s_1, s_2)$} (BA);
	\path [-] (AB) edge[bend right=10] node[pos=0.29] {$(s_2, s_1)$} (BA);
	
	\path [-] (AA) edge[loop, out=170, in=260, min distance=6mm, looseness=50] node[pos=0.5, below=5mm] {$(\id, \id)$} (AA);
	\path [-] (BA) edge[loop, out=280, in=10, min distance=6mm, looseness=50] node[pos=0.5, below=5mm] {$(\id, \id)$} (BA);
	\path [-] (AB) edge[loop, out=100, in=190, min distance=6mm, looseness=50] node[pos=0.5, above=5mm] {$(\id, \id)$} (AB);
	\path [-] (BB) edge[loop, out=350, in=80, min distance=6mm, looseness=50] node[pos=0.5, above=5mm] {$(\id, \id)$} (BB);
\end{scope}
\end{tikzpicture}
}
\caption{The quotient graph $G=\Z^2\backslash\Cay(D_\infty\times D_\infty, S_\ast)$ where $S_\ast=(S\cup\{\id\})^2$ and $S=\{s_1, s_2\}$.}
\label{Fig:A1A1standard}
\end{figure}

Let $\mu$ be a probability measure on $D_\infty$ such that $\supp \mu=S\cup\{\id\}$,
and $\pi^\rho$ be the associated probability measure on $D_\infty \times D_\infty$ for $\rho \in (0, 1]$.
The measure $\pi^\rho$ has the support $S_\ast=(S\cup\{\id\})^2$.
The quotient graph $G=\Z^2\setminus\Cay(D_\infty\times D_\infty, S_\ast)$ is described in Figure \ref{Fig:A1A1standard}.

First we consider the case when $\mu(s_1)$ and $\mu(s_2)$ are equal, i.e.,
\[
\mu(s_1)=\mu(s_2)=\frac{1}{2}(1-\mu(\id)) \quad \text{and} \quad 0<\mu(\id)<1.
\]
In this case, for $v_1=(1, 0), v_2=(0, 1) \in \R^2$, the $1$-forms $(\abr{v_i, \F_\eb})_{\eb \in E(G)}$ for $i=1,2$ are harmonic.
A direct computation yields
\[
\SS^\rho=
\begin{pmatrix}
\frac{1}{4}(1-\mu(\id)) & 0\\
0 & \frac{1}{4}(1-\mu(\id))
\end{pmatrix}
\quad \text{for $\rho \in (0, 1]$}.
\]

Next in the case when $\mu(s_1)$ and $\mu(s_2)$ are not necessarily equal,
the harmonic $1$-forms
\[
u_i(\eb)=\abr{v_i, \F_\eb}-df_i(\eb) \quad \text{for $\eb \in E(G)$ and $i=1, 2$},
\]
are obtained by (possibly non-constant) functions $f_i:V(G) \to \R$.
For $i=1$, let
\[
f_1((0, 0))=f_1((0, 1))=\frac{\mu(s_2)}{2(1-\mu(\id))}
\quad \text{and} \quad
f_1((1, 0))=f_1((1, 1))=\frac{\mu(s_1)}{2(1-\mu(\id))}.
\]
The values of the harmonic $1$-form $u_1$ on the two oriented edges from $(0,0)$ to $(1, 0)$ satisfy the following:
\begin{align*}
u_1((0,0), (0.s_1, 0))+f_1((1,0))-f_1((0,0))&=\F^{(1)}((s_1.o, o))-\F^{(1)}((o,o))=\frac{1}{2},\\
u_1((0,0),(0.s_2, 0))+f_1((1,0))-f_1((0,0))&=\F^{(1)}((s_2.o,o))-\F^{(1)}((o,o))=-\frac{1}{2}.
\end{align*}
Similar identities hold on the four oriented edges from $(0, 0)$ to $(1, 1)$, from $(0, 1)$ to $(1, 0)$, respectively, and on the two oriented edges from $(0, 1)$ to $(1, 1)$.
The values of $u_1$ on the two oriented edges from $(0, 0)$ to $(0, 1)$, and from $(1, 0)$ to $(1, 1)$, respectively, are $0$.
The harmonic $1$-form $u_2$ is obtained analogously.
A direct computation yields
\[
\SS^\rho=
\begin{pmatrix}
\frac{\mu(s_1)\mu(s_2)}{1-\mu(\id)} & 0\\
0 & \frac{\mu(s_1)\mu(s_2)}{1-\mu(\id)}
\end{pmatrix}
\quad \text{for $\rho \in (0, 1]$}.
\]

\subsubsection{Type $\wt A_2$}

Let us consider the affine Weyl group of type $\wt A_2$:
\[
\Gamma=\abr{s_1, s_2, s_3 \mid s_1^2=s_2^2=s_3^2=(s_1s_2)^3=(s_2s_3)^3=(s_3s_1)^3=1}
\]
with $S=\{s_1, s_2, s_3\}$.
The Cayley graph $\Cay(\Gamma, S)$ and the group action on $\R^2$ is described in Figure \ref{Fig:A2} (left).
The group $\Gamma$ is isomorphic to $\Z^2 \rtimes \Sfr_3$ where $\Sfr_3$ is the symmetric group on the set $\{1, 2, 3\}$.
The quotient graph $G=\Z^2\backslash \Cay(\Gamma, S)$ is described in Figure \ref{Fig:A2} (right).

\begin{figure}[h]
\begin{subfigure}{.4\textwidth}
\centering
\scalebox{0.55}[0.55]{
\begin{tikzpicture}[
dot/.style={draw, coordinate},
ver/.style={draw, circle, scale=0.4, fill=black},
]

\node[dot] at (-8,{-4*sqrt(3)}) (-8,{-4*sqrt(3)}) {};
\node[dot] at (-4,{-4*sqrt(3)}) (-4,{-4*sqrt(3)}) {};
\node[dot] at (0,{-4*sqrt(3)}) (0,{-4*sqrt(3)}) {};
\node[dot] at (4,{-4*sqrt(3)}) (4,{-4*sqrt(3)}) {};
\node[dot] at (8,{-4*sqrt(3)}) (8,{-4*sqrt(3)}) {};

\node[dot] at (-6, {-2*sqrt(3)}) (-6, {-2*sqrt(3)}) {};
\node[dot] at (-2, {-2*sqrt(3)}) (-2, {-2*sqrt(3)}) {};
\node[dot] at (2, {-2*sqrt(3)}) (2, {-2*sqrt(3)}) {};
\node[dot] at (6, {-2*sqrt(3)}) (6, {-2*sqrt(3)}) {};
\node[dot] at (10, {-2*sqrt(3)}) (10, {-2*sqrt(3)}) {};

\node[dot] at (-8,0) (-8,0) {};
\node[dot] at (-4,0) (-4,0) {};
\node[dot] at (0,0) (0,0) {};
\node[dot] at (4,0) (4,0) {};
\node[dot] at (8,0) (8,0) {};

\node[dot] at (-6, {2*sqrt(3)}) (-6, {2*sqrt(3)}) {};
\node[dot] at (-2, {2*sqrt(3)}) (-2, {2*sqrt(3)}) {};
\node[dot] at (2, {2*sqrt(3)}) (2, {2*sqrt(3)}) {};
\node[dot] at (6, {2*sqrt(3)}) (6, {2*sqrt(3)}) {};
\node[dot] at (10, {2*sqrt(3)}) (10, {2*sqrt(3)}) {};

\node[dot] at (-8,{4*sqrt(3)}) (-8,{4*sqrt(3)}) {};
\node[dot] at (-4,{4*sqrt(3)}) (-4,{4*sqrt(3)}) {};
\node[dot] at (0,{4*sqrt(3)}) (0,{4*sqrt(3)}) {};
\node[dot] at (4,{4*sqrt(3)}) (4,{4*sqrt(3)}) {};
\node[dot] at (8,{4*sqrt(3)}) (8,{4*sqrt(3)}) {};

\draw[dotted, fill=gray!40] (4,0)--(2,{2*sqrt(3)})--(-2,{2*sqrt(3)})--(-4,0)--(-2,{-2*sqrt(3)})--(2,{-2*sqrt(3)})--cycle;
\draw[dotted, fill=gray!80] (0,0)--(4,0)--(2,{2*sqrt(3)})--cycle;

\draw[dotted] (-8,{-4*sqrt(3)})--(-4,{-4*sqrt(3)});
\draw[dotted] (-4,{-4*sqrt(3)})--(0,{-4*sqrt(3)});
\draw[dotted] (0,{-4*sqrt(3)})--(4,{-4*sqrt(3)});
\draw[dotted] (4,{-4*sqrt(3)})--(8,{-4*sqrt(3)});

\draw[dotted] (-8, {-2*sqrt(3)})--(-6, {-2*sqrt(3)});
\draw[dotted] (-6, {-2*sqrt(3)})--(-2, {-2*sqrt(3)});
\draw[dotted] (-2, {-2*sqrt(3)})--(2, {-2*sqrt(3)});
\draw[dotted] (2, {-2*sqrt(3)})--(6, {-2*sqrt(3)});
\draw[dotted] (6, {-2*sqrt(3)})--(8, {-2*sqrt(3)});

\draw[dotted] (-8,0)--(-4,0);
\draw[dotted] (-4,0)--(0,0);
\draw[dotted] (0,0)--(4,0);
\draw[dotted] (4,0)--(8,0);

\draw[dotted] (-8, {2*sqrt(3)})--(-6, {2*sqrt(3)});
\draw[dotted] (-6, {2*sqrt(3)})--(-2, {2*sqrt(3)});
\draw[dotted] (-2, {2*sqrt(3)})--(2, {2*sqrt(3)});
\draw[dotted] (2, {2*sqrt(3)})--(6, {2*sqrt(3)});
\draw[dotted] (6, {2*sqrt(3)})--(8, {2*sqrt(3)});

\draw[dotted] (-8,{4*sqrt(3)})--(-4,{4*sqrt(3)});
\draw[dotted] (-4,{4*sqrt(3)})--(0,{4*sqrt(3)});
\draw[dotted] (0,{4*sqrt(3)})--(4,{4*sqrt(3)});
\draw[dotted] (4,{4*sqrt(3)})--(8,{4*sqrt(3)});


\draw[dotted] (-8,{-4*sqrt(3)})--(-6, {-2*sqrt(3)});
\draw[dotted] (-4,{-4*sqrt(3)})--(-6, {-2*sqrt(3)});
\draw[dotted] (-4,{-4*sqrt(3)})--(-2, {-2*sqrt(3)});
\draw[dotted] (0,{-4*sqrt(3)})--(-2, {-2*sqrt(3)});
\draw[dotted] (0,{-4*sqrt(3)})--(2, {-2*sqrt(3)});
\draw[dotted] (4,{-4*sqrt(3)})--(2, {-2*sqrt(3)});
\draw[dotted] (4,{-4*sqrt(3)})--(6, {-2*sqrt(3)});
\draw[dotted] (8,{-4*sqrt(3)})--(6, {-2*sqrt(3)});

\draw[dotted] (-8,0)--(-6, {-2*sqrt(3)});
\draw[dotted] (-4,0)--(-6, {-2*sqrt(3)});
\draw[dotted] (-4,0)--(-2, {-2*sqrt(3)});
\draw[dotted] (0,0)--(-2, {-2*sqrt(3)});
\draw[dotted] (0,0)--(2, {-2*sqrt(3)});
\draw[dotted] (4,0)--(2, {-2*sqrt(3)});
\draw[dotted] (4,0)--(6, {-2*sqrt(3)});
\draw[dotted] (8,0)--(6, {-2*sqrt(3)});

\draw[dotted] (-8,0)--(-6, {2*sqrt(3)});
\draw[dotted] (-4,0)--(-6, {2*sqrt(3)});
\draw[dotted] (-4,0)--(-2, {2*sqrt(3)});
\draw[dotted] (0,0)--(-2, {2*sqrt(3)});
\draw[dotted] (0,0)--(2, {2*sqrt(3)});
\draw[dotted] (4,0)--(2, {2*sqrt(3)});
\draw[dotted] (4,0)--(6, {2*sqrt(3)});
\draw[dotted] (8,0)--(6, {2*sqrt(3)});

\draw[dotted] (-8,{4*sqrt(3)})--(-6, {2*sqrt(3)});
\draw[dotted] (-4,{4*sqrt(3)})--(-6, {2*sqrt(3)});
\draw[dotted] (-4,{4*sqrt(3)})--(-2, {2*sqrt(3)});
\draw[dotted] (0,{4*sqrt(3)})--(-2, {2*sqrt(3)});
\draw[dotted] (0,{4*sqrt(3)})--(2, {2*sqrt(3)});
\draw[dotted] (4,{4*sqrt(3)})--(2, {2*sqrt(3)});
\draw[dotted] (4,{4*sqrt(3)})--(6, {2*sqrt(3)});
\draw[dotted] (8,{4*sqrt(3)})--(6, {2*sqrt(3)});


\node[ver] at (-6, {-(10/3)*sqrt(3)}) (-6, {-(10/3)*sqrt(3)}) {};
\node[ver] at (-2, {-(10/3)*sqrt(3)}) (-2, {-(10/3)*sqrt(3)}) {};
\node[ver] at (2, {-(10/3)*sqrt(3)}) (2, {-(10/3)*sqrt(3)}) {};
\node[ver] at (6, {-(10/3)*sqrt(3)}) (6, {-(10/3)*sqrt(3)}) {};

\node[ver] at (-8, {-(8/3)*sqrt(3)}) (-8, {-(8/3)*sqrt(3)}) {};
\node[ver] at (-4, {-(8/3)*sqrt(3)}) (-4, {-(8/3)*sqrt(3)}) {};
\node[ver] at (0, {-(8/3)*sqrt(3)}) (0, {-(8/3)*sqrt(3)}) {};
\node[ver] at (4, {-(8/3)*sqrt(3)}) (4, {-(8/3)*sqrt(3)}) {};
\node[ver] at (8, {-(8/3)*sqrt(3)}) (8, {-(8/3)*sqrt(3)}) {};

\node[ver] at (-6, {-(2/3)*sqrt(3)}) (-6, {-(2/3)*sqrt(3)}) {};
\node[ver] at (-2, {-(2/3)*sqrt(3)}) (-2, {-(2/3)*sqrt(3)}) {};
\node[ver] at (2, {-(2/3)*sqrt(3)}) (2, {-(2/3)*sqrt(3)}) {};
\node[ver] at (6, {-(2/3)*sqrt(3)}) (6, {-(2/3)*sqrt(3)}) {};

\node[ver] at (-8, {-(4/3)*sqrt(3)}) (-8, {-(4/3)*sqrt(3)}) {};
\node[ver] at (-4, {-(4/3)*sqrt(3)}) (-4, {-(4/3)*sqrt(3)}) {};
\node[ver] at (0, {-(4/3)*sqrt(3)}) (0, {-(4/3)*sqrt(3)}) {};
\node[ver] at (4, {-(4/3)*sqrt(3)}) (4, {-(4/3)*sqrt(3)}) {};
\node[ver] at (8, {-(4/3)*sqrt(3)}) (8, {-(4/3)*sqrt(3)}) {};

\node[ver] at (-6, {(2/3)*sqrt(3)}) (-6, {(2/3)*sqrt(3)}) {};
\node[ver] at (-2, {(2/3)*sqrt(3)}) (-2, {(2/3)*sqrt(3)}) {};
\node[ver] at (2, {(2/3)*sqrt(3)}) (2, {(2/3)*sqrt(3)}) {}; 
\node[] at ({2+0.3}, {(2/3)*sqrt(3)-0.2}) ({2+0.3}, {(2/3)*sqrt(3)-0.2}) {$o$}; 
\node[ver] at (6, {(2/3)*sqrt(3)}) (6, {(2/3)*sqrt(3)}) {};

\node[ver] at (-8, {(4/3)*sqrt(3)}) (-8, {(4/3)*sqrt(3)}) {};
\node[ver] at (-4, {(4/3)*sqrt(3)}) (-4, {(4/3)*sqrt(3)}) {};
\node[ver] at (0, {(4/3)*sqrt(3)}) (0, {(4/3)*sqrt(3)}) {};
\node[ver] at (4, {(4/3)*sqrt(3)}) (4, {(4/3)*sqrt(3)}) {};
\node[ver] at (8, {(4/3)*sqrt(3)}) (8, {(4/3)*sqrt(3)}) {};

\node[ver] at (-6, {(10/3)*sqrt(3)}) (-6, {(10/3)*sqrt(3)}) {};
\node[ver] at (-2, {(10/3)*sqrt(3)}) (-2, {(10/3)*sqrt(3)}) {};
\node[ver] at (2, {(10/3)*sqrt(3)}) (2, {(10/3)*sqrt(3)}) {};
\node[ver] at (6, {(10/3)*sqrt(3)}) (6, {(10/3)*sqrt(3)}) {};

\node[ver] at (-8, {(8/3)*sqrt(3)}) (-8, {(8/3)*sqrt(3)}) {};
\node[ver] at (-4, {(8/3)*sqrt(3)}) (-4, {(8/3)*sqrt(3)}) {};
\node[ver] at (0, {(8/3)*sqrt(3)}) (0, {(8/3)*sqrt(3)}) {};
\node[ver] at (4, {(8/3)*sqrt(3)}) (4, {(8/3)*sqrt(3)}) {};
\node[ver] at (8, {(8/3)*sqrt(3)}) (8, {(8/3)*sqrt(3)}) {};



\draw[] (-8, {(8/3)*sqrt(3)})--(-6, {(10/3)*sqrt(3)});
\draw[] (-6, {(10/3)*sqrt(3)})--(-4, {(8/3)*sqrt(3)});
\draw[] (-4, {(8/3)*sqrt(3)})--(-2, {(10/3)*sqrt(3)});
\draw[] (-2, {(10/3)*sqrt(3)})--(0, {(8/3)*sqrt(3)});
\draw[] (0, {(8/3)*sqrt(3)})--(2, {(10/3)*sqrt(3)});
\draw[] (2, {(10/3)*sqrt(3)})--(4, {(8/3)*sqrt(3)});
\draw[] (4, {(8/3)*sqrt(3)})--(6, {(10/3)*sqrt(3)});
\draw[] (6, {(10/3)*sqrt(3)})--(8, {(8/3)*sqrt(3)});

\draw[] (-8, {(4/3)*sqrt(3)})--(-6, {(2/3)*sqrt(3)});
\draw[] (-6, {(2/3)*sqrt(3)})--(-4, {(4/3)*sqrt(3)});
\draw[] (-4, {(4/3)*sqrt(3)})--(-2, {(2/3)*sqrt(3)});
\draw[] (-2, {(2/3)*sqrt(3)})--(0, {(4/3)*sqrt(3)});
\draw[] (0, {(4/3)*sqrt(3)})--(2, {(2/3)*sqrt(3)}) node[pos=0.4, above=0.1mm] {$s_2$};
\draw[] (2, {(2/3)*sqrt(3)})--(4, {(4/3)*sqrt(3)}) node[pos=0.6, above=0.1mm] {$s_3$};
\draw[] (4, {(4/3)*sqrt(3)})--(6, {(2/3)*sqrt(3)});
\draw[] (6, {(2/3)*sqrt(3)})--(8, {(4/3)*sqrt(3)});

\draw[] (-8, {-(4/3)*sqrt(3)})--(-6, {-(2/3)*sqrt(3)});
\draw[] (-6, {-(2/3)*sqrt(3)})--(-4, {-(4/3)*sqrt(3)});
\draw[] (-4, {-(4/3)*sqrt(3)})--(-2, {-(2/3)*sqrt(3)});
\draw[] (-2, {-(2/3)*sqrt(3)})--(0, {-(4/3)*sqrt(3)});
\draw[] (0, {-(4/3)*sqrt(3)})--(2, {-(2/3)*sqrt(3)});
\draw[] (2, {-(2/3)*sqrt(3)})--(4, {-(4/3)*sqrt(3)});
\draw[] (4, {-(4/3)*sqrt(3)})--(6, {-(2/3)*sqrt(3)});
\draw[] (6, {-(2/3)*sqrt(3)})--(8, {-(4/3)*sqrt(3)});

\draw[] (-8, {-(8/3)*sqrt(3)})--(-6, {-(10/3)*sqrt(3)});
\draw[] (-6, {-(10/3)*sqrt(3)})--(-4, {-(8/3)*sqrt(3)});
\draw[] (-4, {-(8/3)*sqrt(3)})--(-2, {-(10/3)*sqrt(3)});
\draw[] (-2, {-(10/3)*sqrt(3)})--(0, {-(8/3)*sqrt(3)});
\draw[] (0, {-(8/3)*sqrt(3)})--(2, {-(10/3)*sqrt(3)});
\draw[] (2, {-(10/3)*sqrt(3)})--(4, {-(8/3)*sqrt(3)});
\draw[] (4, {-(8/3)*sqrt(3)})--(6, {-(10/3)*sqrt(3)});
\draw[] (6, {-(10/3)*sqrt(3)})--(8, {-(8/3)*sqrt(3)});


\draw[] (-6, {(10/3)*sqrt(3)})--(-6, {(12/3)*sqrt(3)});
\draw[] (-2, {(10/3)*sqrt(3)})--(-2, {(12/3)*sqrt(3)});
\draw[] (2, {(10/3)*sqrt(3)})--(2, {(12/3)*sqrt(3)});
\draw[] (6, {(10/3)*sqrt(3)})--(6, {(12/3)*sqrt(3)});

\draw[] (-8, {(4/3)*sqrt(3)})--(-8, {(8/3)*sqrt(3)});
\draw[] (-4, {(4/3)*sqrt(3)})--(-4, {(8/3)*sqrt(3)});
\draw[] (0, {(4/3)*sqrt(3)})--(0, {(8/3)*sqrt(3)});
\draw[] (4, {(4/3)*sqrt(3)})--(4, {(8/3)*sqrt(3)});
\draw[] (8, {(4/3)*sqrt(3)})--(8, {(8/3)*sqrt(3)});

\draw[] (-6, {(2/3)*sqrt(3)})--(-6, {-(2/3)*sqrt(3)});
\draw[] (-2, {(2/3)*sqrt(3)})--(-2, {-(2/3)*sqrt(3)});
\draw[] (2, {(2/3)*sqrt(3)})--(2, {-(2/3)*sqrt(3)}) node[pos=0.7, right=0.1mm] {$s_1$};
\draw[] (6, {(2/3)*sqrt(3)})--(6, {-(2/3)*sqrt(3)});

\draw[] (-8, {-(4/3)*sqrt(3)})--(-8, {-(6/3)*sqrt(3)});
\draw[] (-4, {-(4/3)*sqrt(3)})--(-4, {-(6/3)*sqrt(3)});
\draw[] (0, {-(4/3)*sqrt(3)})--(0, {-(6/3)*sqrt(3)});
\draw[] (4, {-(4/3)*sqrt(3)})--(4, {-(6/3)*sqrt(3)});
\draw[] (8, {-(4/3)*sqrt(3)})--(8, {-(6/3)*sqrt(3)});

\draw[] (-8, {-(4/3)*sqrt(3)})--(-8, {-(8/3)*sqrt(3)});
\draw[] (-4, {-(4/3)*sqrt(3)})--(-4, {-(8/3)*sqrt(3)});
\draw[] (0, {-(4/3)*sqrt(3)})--(0, {-(8/3)*sqrt(3)});
\draw[] (4, {-(4/3)*sqrt(3)})--(4, {-(8/3)*sqrt(3)});
\draw[] (8, {-(4/3)*sqrt(3)})--(8, {-(8/3)*sqrt(3)});

\draw[] (-6, {-(10/3)*sqrt(3)})--(-6, {-(12/3)*sqrt(3)});
\draw[] (-2, {-(10/3)*sqrt(3)})--(-2, {-(12/3)*sqrt(3)});
\draw[] (2, {-(10/3)*sqrt(3)})--(2, {-(12/3)*sqrt(3)});
\draw[] (6, {-(10/3)*sqrt(3)})--(6, {-(12/3)*sqrt(3)});

\end{tikzpicture}
}
\end{subfigure}
\qquad \qquad \qquad
\begin{subfigure}{.4\textwidth}
\centering
\scalebox{0.55}[0.55]{
\begin{tikzpicture}[
dot/.style={draw, coordinate},
ver/.style={draw, circle, scale=0.4, fill=black},
every loop/.style={}
]

\node[ver] at (-2, {(2/3)*sqrt(3)}) (-2, {(2/3)*sqrt(3)}) {};
\node[ver] at (-2, {-(2/3)*sqrt(3)}) (-2, {-(2/3)*sqrt(3)}) {};

\node[ver] at (2, {(2/3)*sqrt(3)}) (2, {(2/3)*sqrt(3)}) {}; 
\node[ver] at (2, {-(2/3)*sqrt(3)}) (2, {-(2/3)*sqrt(3)}) {};

\node[ver] at (0, {(4/3)*sqrt(3)}) (0, {(4/3)*sqrt(3)}) {};
\node[ver] at (0, {-(4/3)*sqrt(3)}) (0, {-(4/3)*sqrt(3)}) {};

\draw[] (2, {(2/3)*sqrt(3)})--(2, {-(2/3)*sqrt(3)}) node[pos=0.5, right=0.1mm] {$s_1$};
\draw[] (-2, {(2/3)*sqrt(3)})--(-2, {-(2/3)*sqrt(3)}) node[pos=0.5, left=0.1mm] {$s_2$};
\draw[] (-2, {(2/3)*sqrt(3)})--(0, {(4/3)*sqrt(3)}) node[pos=0.5, above=0.1mm] {$s_1$};
\draw[] (2, {(2/3)*sqrt(3)})--(0, {(4/3)*sqrt(3)}) node[pos=0.5, above=0.1mm] {$s_2$};
\draw[] (-2, {-(2/3)*sqrt(3)})--(0, {-(4/3)*sqrt(3)}) node[pos=0.5, below=0.1mm] {$s_1$};
\draw[] (2, {-(2/3)*sqrt(3)})--(0, {-(4/3)*sqrt(3)}) node[pos=0.5, below=0.1mm] {$s_2$};

\draw[] (2, {(2/3)*sqrt(3)})--(-2, {-(2/3)*sqrt(3)}) node[pos=0.3, above=0.1mm] {$s_3$};
\draw[] (2, {-(2/3)*sqrt(3)})--(-2, {(2/3)*sqrt(3)}) node[pos=0.3, below=0.1mm] {$s_3$};
\draw[] (0, {(4/3)*sqrt(3)})--(0, {-(4/3)*sqrt(3)}) node[pos=0.25, left=0.1mm] {$s_3$};

\end{tikzpicture}
}
\end{subfigure}
\caption{The Cayley graph of the affine Weyl group $(\Gamma, S)$ of type $\wt A_2$ and the action on $\R^2$ (left), 
and the quotient graph $G=\Z^2\backslash \Cay(\Gamma, S)$ (right).}
\label{Fig:A2}
\end{figure}

Let us consider $\F:\Gamma \to \R^2$ such that the origin is the barycenter of a fundamental chamber (which is an equilateral triangle of side length $(2\sqrt{3})^{1/2}/3$).
A fundamental domain for the $\Z^2$-action is a hexagon of unit area.
Let $\mu$ be a probability measure on $\Gamma$ such that
\[
\mu(s_1)=\mu(s_2)=\mu(s_3)=\frac{1}{3}(1-\mu(\id)) \quad \text{and} \quad 0<\mu(\id)<1.
\]
The matrix for $\mu$ in the local central limit theorem (Theorem \ref{Thm:LCLT}) is computed as
\[
\SS^\mu=
\begin{pmatrix}
\frac{\sqrt{3}}{27}(1-\mu(\id)) & 0\\
0 & \frac{\sqrt{3}}{27}(1-\mu(\id))
\end{pmatrix}.
\]

\subsubsection{Type $\wt C_2$}

Let us consider the affine Weyl group of type $\wt C_2$:
\[
\Gamma=\abr{s_1, s_2, s_3 \mid s_1^2=s_2^2=s_3^2=(s_1s_2)^4=(s_2s_3)^4=(s_3s_1)^2=1}
\]
with $S=\{s_1, s_2, s_3\}$.
The Cayley graph $\Cay(\Gamma, S)$ and the group action on $\R^2$ is described in Figure \ref{Fig:C2} (left).
The group $\Gamma$ is isomorphic to $\Z^2 \rtimes ((\Z/2)^2\rtimes\Sfr_2)$ where $(\Z/2)^2\rtimes\Sfr_2$ is the signed permutations on the set $\{1, 2\}$.
The quotient graph $G$ of the Cayley graph with the set of generators $S$ by the lattice $\Z^2$ is described in Figure \ref{Fig:C2} (right).

\begin{figure}[h]
\begin{subfigure}{.4\textwidth}
\centering
\scalebox{0.55}[0.55]{
\begin{tikzpicture}[
dot/.style={draw, coordinate},
ver/.style={draw, circle, scale=0.4, fill=black},
]

\draw[dotted, fill=gray!40] (3,3)--(-3,3)--(-3,-3)--(3,-3)--cycle;
\draw[dotted, fill=gray!80] (0,0)--(3,0)--(3,3)--cycle;

\node[dot] at (-6,6) (-6,6) {};
\node[dot] at (-3,6) (-3,6) {};
\node[dot] at (0,6) (0,6) {};
\node[dot] at (3,6) (3,6) {};
\node[dot] at (6,6) (6,6) {};

\node[dot] at (-6,3) (-6,3) {};
\node[dot] at (-3,3) (-3,3) {};
\node[dot] at (0,3) (0,3) {};
\node[dot] at (3,3) (3,3) {};
\node[dot] at (6,3) (6,3) {};

\node[dot] at (-6,0) (-6,0) {};
\node[dot] at (-3,0) (-3,0) {};
\node[dot] at (0,0) (0,0) {};
\node[dot] at (3,0) (3,0) {};
\node[dot] at (6,0) (6,0) {};

\node[dot] at (-6,-3) (-6,-3) {};
\node[dot] at (-3,-3) (-3,-3) {};
\node[dot] at (0,-3) (0,-3) {};
\node[dot] at (3,-3) (3,-3) {};
\node[dot] at (6,-3) (6,-3) {};

\node[dot] at (-6,-6) (-6,-6) {};
\node[dot] at (-3,-6) (-3,-6) {};
\node[dot] at (0,-6) (0,-6) {};
\node[dot] at (3,-6) (3,-6) {};
\node[dot] at (6,-6) (6,-6) {};


\draw[dotted] (-6,6) -- (-3,6);
\draw[dotted] (-3,6) -- (0,6);
\draw[dotted] (0,6) -- (3,6);
\draw[dotted] (3,6) -- (6,6);

\draw[dotted] (-6,3) -- (-3,3);
\draw[dotted] (-3,3) -- (0,3);
\draw[dotted] (0,3) -- (3,3);
\draw[dotted] (3,3) -- (6,3);

\draw[dotted] (-6,0) -- (-3,0);
\draw[dotted] (-3,0) -- (0,0);
\draw[dotted] (0,0) -- (3,0);
\draw[dotted] (3,0) -- (6,0);

\draw[dotted] (-6,-3) -- (-3,-3);
\draw[dotted] (-3,-3) -- (0,-3);
\draw[dotted] (0,-3) -- (3,-3);
\draw[dotted] (3,-3) -- (6,-3);

\draw[dotted] (-6,-6) -- (-3,-6);
\draw[dotted] (-3,-6) -- (0,-6);
\draw[dotted] (0,-6) -- (3,-6);
\draw[dotted] (3,-6) -- (6,-6);


\draw[dotted] (-6,6) -- (-6,3);
\draw[dotted] (-3,6) -- (-3,3);
\draw[dotted] (0,6) -- (0,3);
\draw[dotted] (3,6) -- (3,3);
\draw[dotted] (6,6) -- (6,3);

\draw[dotted] (-6,0) -- (-6,3);
\draw[dotted] (-3,0) -- (-3,3);
\draw[dotted] (0,0) -- (0,3);
\draw[dotted] (3,0) -- (3,3);
\draw[dotted] (6,0) -- (6,3);

\draw[dotted] (-6,0) -- (-6,-3);
\draw[dotted] (-3,0) -- (-3,-3);
\draw[dotted] (0,0) -- (0,-3);
\draw[dotted] (3,0) -- (3,-3);
\draw[dotted] (6,0) -- (6,-3);

\draw[dotted] (-6,-6) -- (-6,-3);
\draw[dotted] (-3,-6) -- (-3,-3);
\draw[dotted] (0,-6) -- (0,-3);
\draw[dotted] (3,-6) -- (3,-3);
\draw[dotted] (6,-6) -- (6,-3);


\draw[dotted] (-6,6) -- (-3,3);
\draw[dotted] (0,6) -- (-3,3);
\draw[dotted] (0,6) -- (3,3);
\draw[dotted] (6,6) -- (3,3);

\draw[dotted] (-6,0) -- (-3,3);
\draw[dotted] (0,0) -- (-3,3);
\draw[dotted] (0,0) -- (3,3);
\draw[dotted] (6,0) -- (3,3);

\draw[dotted] (-6,0) -- (-3,-3);
\draw[dotted] (0,0) -- (-3,-3);
\draw[dotted] (0,0) -- (3,-3);
\draw[dotted] (6,0) -- (3,-3);

\draw[dotted] (-6,-6) -- (-3,-3);
\draw[dotted] (0,-6) -- (-3,-3);
\draw[dotted] (0,-6) -- (3,-3);
\draw[dotted] (6,-6) -- (3,-3);



\node[ver] at (-4,5) (-4,5) {};
\node[ver] at (-2,5) (-2,5) {};
\node[ver] at (2,5) (2,5) {};
\node[ver] at (4,5) (4,5) {};

\node[ver] at (-5,4) (-5,4) {};
\node[ver] at (-1,4) (-1,4) {};
\node[ver] at (1,4) (1,4) {};
\node[ver] at (5,4) (5,4) {};

\node[ver] at (-5,2) (-5,2) {};
\node[ver] at (-1,2) (-1,2) {};
\node[ver] at (1,2) (1,2) {};
\node[ver] at (5,2) (5,2) {};

\node[ver] at (-4,1) (-4,1) {};
\node[ver] at (-2,1) (-2,1) {};
\node[ver] at (2,1) (2,1) {}; 
\node[above] at (2,1) {$o$}; 
\node[ver] at (4,1) (4,1) {};

\node[ver] at (-4,-1) (-4,-1) {};
\node[ver] at (-2,-1) (-2,-1) {};
\node[ver] at (2,-1) (2,-1) {};
\node[ver] at (4,-1) (4,-1) {};

\node[ver] at (-5,-2) (-5,-2) {};
\node[ver] at (-1,-2) (-1,-2) {};
\node[ver] at (1,-2) (1,-2) {};
\node[ver] at (5,-2) (5,-2) {};

\node[ver] at (-5,-4) (-5,-4) {};
\node[ver] at (-1,-4) (-1,-4) {};
\node[ver] at (1,-4) (1,-4) {};
\node[ver] at (5,-4) (5,-4) {};

\node[ver] at (-4,-5) (-4,-5) {};
\node[ver] at (-2,-5) (-2,-5) {};
\node[ver] at (2,-5) (2,-5) {};
\node[ver] at (4,-5) (4,-5) {};

\draw[] (-4,5)--(-2,5);
\draw[] (2,5)--(4,5);

\draw[] (-1,2)--(1,2);
\draw[] (-1,4)--(1,4);
\draw[] (-5,2)--(-6,2);
\draw[] (-5,4)--(-6,4);
\draw[] (5,2)--(6,2);
\draw[] (5,4)--(6,4);

\draw[] (-1,-2)--(1,-2);
\draw[] (-1,-4)--(1,-4);
\draw[] (-5,-2)--(-6,-2);
\draw[] (-5,-4)--(-6,-4);
\draw[] (5,-2)--(6,-2);
\draw[] (5,-4)--(6,-4);

\draw[] (-4,1)--(-2,1);
\draw[] (2,1)--(4,1) node[pos=0.7, above=0.1mm] {$s_3$};
\draw[] (-4,-1)--(-2,-1);
\draw[] (2,-1)--(4,-1);

\draw[] (-4,-5)--(-2,-5);
\draw[] (2,-5)--(4,-5);

\draw[] (1,2)--(1,4);
\draw[] (-1,2)--(-1,4);
\draw[] (5,2)--(5,4);
\draw[] (-5,2)--(-5,4);

\draw[] (1,-2)--(1,-4);
\draw[] (-1,-2)--(-1,-4);
\draw[] (5,-2)--(5,-4);
\draw[] (-5,-2)--(-5,-4);

\draw[] (2,1)--(2,-1) node[pos=0.7, right=0.1mm] {$s_1$};
\draw[] (4,1)--(4,-1);
\draw[] (-2,1)--(-2,-1);
\draw[] (-4,1)--(-4,-1);

\draw[] (2,5)--(2,6);
\draw[] (4,5)--(4,6);
\draw[] (-2,5)--(-2,6);
\draw[] (-4,5)--(-4,6);

\draw[] (2,-5)--(2,-6);
\draw[] (4,-5)--(4,-6);
\draw[] (-2,-5)--(-2,-6);
\draw[] (-4,-5)--(-4,-6);


\draw[] (-4,5)--(-5,4);
\draw[] (-2,5)--(-1,4);
\draw[] (2,5)--(1,4);
\draw[] (4,5)--(5,4);

\draw[] (-4,1)--(-5,2);
\draw[] (-2,1)--(-1,2);
\draw[] (2,1)--(1,2) node[pos=0.6, above=0.1mm] {$s_2$};
\draw[] (4,1)--(5,2);

\draw[] (-4,-1)--(-5,-2);
\draw[] (-2,-1)--(-1,-2);
\draw[] (2,-1)--(1,-2);
\draw[] (4,-1)--(5,-2);

\draw[] (-4,-5)--(-5,-4);
\draw[] (-2,-5)--(-1,-4);
\draw[] (2,-5)--(1,-4);
\draw[] (4,-5)--(5,-4);

\end{tikzpicture}
}
\end{subfigure}
\qquad \qquad
\begin{subfigure}{.4\textwidth}
\centering
\scalebox{0.55}[0.55]{
\begin{tikzpicture}[
dot/.style={draw, coordinate},
ver/.style={draw, circle, scale=0.4, fill=black},
every loop/.style={}
]

\node[ver] at (-2,1) (-2,1) {};
\node[ver] at (-1,2) (-1,2) {};
\node[ver] at (1,2) (1,2) {};
\node[ver] at (2,1) (2,1) {};

\node[ver] at (-2,-1) (-2,-1) {};
\node[ver] at (-1,-2) (-1,-2) {};
\node[ver] at (1,-2) (1,-2) {};
\node[ver] at (2,-1) (2,-1) {};


\draw[] (1,2)--(-1,2) node[pos=0.5, above=0.1mm] {$s_1$};
\draw[] (1,-2)--(-1,-2) node[pos=0.5, below=0.1mm] {$s_1$};
\draw[] (2,1)--(2,-1) node[pos=0.5, right=0.1mm] {$s_1$};
\draw[] (-2,1)--(-2,-1) node[pos=0.5, left=0.1mm] {$s_1$};

\draw[] (1,2)--(2,1) node[pos=0.6, above=0.1mm] {$s_2$};
\draw[] (-1,2)--(-2,1) node[pos=0.6, above=0.1mm] {$s_2$};
\draw[] (1,-2)--(2,-1) node[pos=0.6, below=0.1mm] {$s_2$};
\draw[] (-1,-2)--(-2,-1) node[pos=0.6, below=0.1mm] {$s_2$};

\draw[] (2,1)--(-2,1) node[pos=0.5, above=0.1mm] {$s_3$};
\draw[] (2,-1)--(-2,-1) node[pos=0.5, below=0.1mm] {$s_3$};
\draw[] (1,2)--(1,-2) node[pos=0.5, right=0.1mm] {$s_3$};
\draw[] (-1,2)--(-1,-2) node[pos=0.5, left=0.1mm] {$s_3$};

\end{tikzpicture}
}
\end{subfigure}
\caption{The Cayley graph of the affine Weyl group $(\Gamma, S)$ of type $\wt C_2$ and the action on $\R^2$ (left),
and the quotient graph $G=\Z^2\backslash \Cay(\Gamma, S)$ (right).}
\label{Fig:C2}
\end{figure}

Let us consider $\F:\Gamma \to \R^2$ such that the origin is the barycenter of a square of side length $1/4$ (where a fundamental chamber is an isosceles right triangle of equal side length $1/2$).
A fundamental domain for the $\Z^2$-action has unit area.
Let $\mu$ be a probability measure on $\Gamma$ such that
\[
\mu(s_1)=\mu(s_2)=\mu(s_3)=\frac{1}{3}(1-\mu(\id)) \quad \text{and} \quad 0<\mu(\id)<1.
\]
The matrix for $\mu$ in the local central limit theorem (Theorem \ref{Thm:LCLT}) is computed as
\[
\SS^\mu=
\begin{pmatrix}
\frac{1}{24}(1-\mu(\id)) & 0\\
0 & \frac{1}{24}(1-\mu(\id))
\end{pmatrix}.
\]

\appendix

\section{Noise sensitivity problem on $\Z^m$}\label{Sec:appendix}

\begin{theorem}\label{Thm:Z}
Let $m$ be a positive integer and $\mu$ be a probability measure on $\Z^m$.
If $\mu$ has a finite second moment and the support generates the group as a semigroup,
then
\[
\lim_{\rho \to 1}\limsup_{n \to \infty}\|\pi^\rho_n-\mu_n\times \mu_n\|_\TV=0
\quad
{and}
\quad
\lim_{\rho \to 0}\liminf_{n \to \infty}\|\pi^\rho_n-\mu_n\times \mu_n\|_\TV =1.
\]
\end{theorem}

Let $\SS$ be a non-degenerate positive definite (covariance) matrix of size $m$ and
\[
f_\SS(v):=\frac{1}{\sqrt{(2\pi)^m \det \SS}}e^{-\frac{1}{2}\abr{v, \SS^{-1}v}} \quad \text{for $v \in \R^m$}.
\]
Let us define the function
$F(v):=\sum_{x\in \Z^m}f_\SS(x+v)$ for $v \in \R^m$.
Note that since $F$ is $\Z^m$-periodic on $\R^m$, it is regarded as a function on $[0, 1)^m$.
Let
\[
\wh F(x):=\int_{[0, 1)^m}F(v)e^{-2\pi i \abr{v, x}}\,dv
\quad \text{and} \quad
\wh f_\SS(x):=\int_{\R^m}f_\SS(v)e^{-2\pi i \abr{v, x}}\,dv \quad \text{for $x \in \Z^m$}.
\]
Since $f_\SS$ is in the Schwartz class on $\R^m$, we have $\wh F(x)=\wh f_\SS(x)$ for $x \in \Z^m$ and
\[
F(v)=\sum_{x \in \Z^m}\wh f_\SS(x)e^{2\pi i \abr{x, v}},
\]
where the right hand side is absolutely convergent.
In the case when $v=0$, the Poisson summation formula is obtained by a direct computation,
\[
\sum_{x \in \Z^m}f_\SS(x)=\sum_{x \in \Z^m}\wh f_\SS(x)=\sum_{x \in \Z^m}e^{-2\pi^2\abr{x, \SS x}}.
\]
Thus there exists a constant $c_\SS>0$ such that for all $n \in \Z_{>0}$,
\begin{equation}\label{Eq:Poisson}
\sum_{x \in \Z^m}f_{n\SS}(x)=\sum_{x \in \Z^m}e^{-2\pi^2n\abr{x, \SS x}}=1+O_\SS(e^{-c_\SS n}).
\end{equation}

Let $\SS_1$ and $\SS_2$ be covariance matrices of size $m$.
First we consider the upper bound.
For every real $\lambda>0$ and for every $n \in \Z_{>0}$,
\begin{align}\label{Eq:sum}
&\sum_{\|x\| \le \lambda n^{1/2}}|f_{n\SS_1}(x)-f_{n\SS_2}(x)| \nonumber\\
&\le\sum_{\|x\| \le \lambda n^{1/2}}\left|f_{n\SS_1}(x)-\sqrt{\frac{\det \SS_2}{\det \SS_1}}f_{n\SS_2}(x)\right|+\sum_{\|x\| \le \lambda n^{1/2}}\left|1-\sqrt{\frac{\det \SS_2}{\det \SS_1}}\right|f_{n\SS_2}(x).
\end{align}
The first sum in \eqref{Eq:sum} is estimated as follows:
For $x \in \Z^m$,
\begin{align*}
\Big|1-e^{-\frac{1}{2n}\abr{x, (\SS_2^{-1}-\SS_1^{-1})x}}\Big|
\le e^{\frac{1}{2n}|\abr{x, (\SS_2^{-1}-\SS_1^{-1})x}|}- e^{-\frac{1}{2n}|\abr{x, (\SS_2^{-1}-\SS_1^{-1})x}|}.
\end{align*}
Since $|\abr{x, (\SS_2^{-1}-\SS_1^{-1})x}|\le \|\SS_2^{-1}-\SS_1^{-1}\|\|x\|^2$,
for $x \in \Z^m$ with $\|x\|\le \lambda n^{1/2}$,
\[
\Big|1-e^{-\frac{1}{2n}\abr{x, (\SS_2^{-1}-\SS_1^{-1})x}}\Big| \le e^{\frac{\lambda^2}{2}\|\SS_2^{-1}-\SS_1^{-1}\|}-e^{-\frac{\lambda^2}{2}\|\SS_2^{-1}-\SS_1^{-1}\|}=2\sinh\(\frac{\lambda^2}{2}\|\SS_2^{-1}-\SS_1^{-1}\|\).
\]
Therefore the first sum in the right hand side of \eqref{Eq:sum} is at most
\begin{align*}
\sum_{\|x\| \le \lambda n^{1/2}}
\Big|1-e^{-\frac{1}{2n}\abr{x, (\SS_2^{-1}-\SS_1^{-1})x}}\Big|f_{n\SS_1}(x) \le 3\lambda^2\|\SS_2^{-1}-\SS_1^{-1}\|\(1+O_{\SS_1}(e^{-c_{\SS_1}n})\),
\end{align*} 
if $\|\SS_2^{-1}-\SS_1^{-1}\| \le 2/\lambda^2$ for all $n \in \Z_{>0}$ and for all $\lambda>0$ by \eqref{Eq:Poisson}.
The second sum in the right hand side of \eqref{Eq:sum} is at most, for all $n \in \Z_{>0}$, by \eqref{Eq:Poisson}, 
\begin{align*}
\left|1-\sqrt{\frac{\det \SS_2}{\det \SS_1}}\right|\sum_{x \in \Z^m}
f_{n\SS_2}(x)
=\left|1-\sqrt{\frac{\det \SS_2}{\det \SS_1}}\right|\(1+O_{\SS_2}(e^{-c_{\SS_2}n})\).
\end{align*}
Summarizing the above estimates yields for every $\lambda>0$,
there exist constants $C_{1, 2, \lambda}=C_{\SS_1, \SS_2, \lambda}>0$ and $c_{1, 2}=c_{\SS_1, \SS_2}>0$ such that for all $n \in \Z_{>0}$,
\begin{align}\label{Eq:limsup_rho}
\sum_{\|x\| \le \lambda n^{1/2}}|f_{n\SS_1}(x)-f_{n\SS_2}(x)|
\le 3\lambda^2\|\SS_2^{-1}-\SS_1^{-1}\|+\left|1-\sqrt{\frac{\det \SS_2}{\det \SS_1}}\right|+C_{1, 2, \lambda} e^{-c_{1,2}n}.
\end{align}

Next let us consider the lower bound.
Noting that $\sqrt{f_1-f_2} \ge \sqrt{f_1}-\sqrt{f_2}$ for $f_1 \ge f_2 \ge 0$,
we have by squaring both sides and summing over $\Z^m$,
\begin{align*}
\|f_{n\SS_1}-f_{n\SS_2}\|_1 \ge \sum_{x \in \Z^m}f_{n\SS_1}(x)+\sum_{x\in \Z^m}f_{n\SS_2}(x)-2\sum_{x \in \Z^m}\sqrt{f_{n\SS_1}(x)f_{n\SS_2}(x)}.
\end{align*}
In the following, we assume that $\SS_1^{-1}+\SS_2^{-1}$ is invertible.
In that case, by \eqref{Eq:Poisson},
\begin{align*}
\sum_{x \in \Z^m}\sqrt{f_{n\SS_1}(x)f_{n\SS_2}(x)}
&=\sum_{x \in \Z^m}\frac{1}{(2\pi)^{\frac{m}{2}}(\det(n\SS_1)\det(n\SS_2))^{\frac{1}{4}}}
e^{-\frac{1}{4n}\abr{x, (\SS_1^{-1}+\SS_2^{-1})x}}\\
&=\frac{2^{\frac{m}{2}}(\det\((\SS_1^{-1}+\SS_2^{-1})^{-1}\))^{\frac{1}{2}}}{(\det \SS_1 \det \SS_2)^{\frac{1}{4}}}
\sum_{x \in \Z^m}e^{-4\pi^2 n\abr{x, (\SS_1^{-1}+\SS_2^{-1})^{-1}x}}.
\end{align*}
Let us focus on the special case when $m=2$, and
\[
\SS^1=\sigma^2
\begin{pmatrix}
1 & 0\\
0 & 1
\end{pmatrix}
\quad \text{and} \quad 
\SS^\rho=\sigma^2
\begin{pmatrix}
1 & 1-\rho\\
1-\rho & 1
\end{pmatrix}
\quad \text{for $\sigma>0$ and $0<\rho \le 1$}.
\]
In this case,
\begin{align*}
(\SS^1)^{-1}+(\SS^\rho)^{-1}
=\frac{1}{\sigma^2}
\begin{pmatrix}
1+\frac{1}{\rho(2-\rho)} & -\frac{1-\rho}{\rho(2-\rho)}\\
-\frac{1-\rho}{\rho(2-\rho)} & 1+\frac{1}{\rho(2-\rho)}
\end{pmatrix}.
\end{align*}
Furthermore, 
$\det \SS^1=\sigma^4$, $\det \SS^\rho=\sigma^4\rho(2-\rho)$ and
\[
\det\((\SS^1)^{-1}+(\SS^\rho)^{-1}\)=\frac{1}{\sigma^4}\(\(1+\frac{1}{\rho(2-\rho)}\)^2-\(\frac{1-\rho}{\rho(2-\rho)}\)^2\).
\]
Hence one computes for all small enough $\rho>0$,
\begin{align*}
\frac{2}{(\rho(2-\rho))^{\frac{1}{4}}\(\(1+\frac{1}{\rho(2-\rho)}\)^2-\(\frac{1-\rho}{\rho(2-\rho)}\)^2\)^{\frac{1}{2}}}
=\frac{2(\rho(2-\rho))^{\frac{3}{4}}}{\sqrt{(\rho(2-\rho)+1)^2-(1-\rho)^2}}<2\rho^{\frac{1}{4}}.
\end{align*}
Summarizing the above computations and letting $n \to \infty$ yield, for all small enough $\rho>0$,
\begin{equation}\label{Eq:liminf_rho}
\liminf_{n \to \infty}\|f_{n\SS^\rho}-f_{n\SS^1}\|_1 \ge 2-4\rho^{\frac{1}{4}}.
\end{equation}

\proof[Proof of Theorem \ref{Thm:Z}]
Noting that $\E_{\pi^\rho_n} \wb_n$ is independent of $\rho \in [0, 1]$,
we assume that $\mu$ has mean zero up to shifting by the mean.
Further, we assume that $\mu$ is aperiodic (i.e., for each $x\in \Z^m$ for all large enough $n$ one has $\mu_n(x)>0$) and the general case is reduced to this case by dividing the times according to the period.
For $0<\rho\le 1$, we note that $\pi^\rho$ is aperiodic since $\pi^\rho$ and $\mu \times \mu$ have the same support in $\Z^{2m}$.
Let $\SS^\rho$ and $\SS^1$ denote the covariance matrices of size $2m$ for $\pi^\rho$ and $\mu\times \mu$ respectively.

First let us show the upper bound in the claim.
For all real $\lambda>0$ and integers $n \in \Z_{>0}$, 
\begin{align*}
\sum_{\|\xb\|>\lambda n^{1/2}}\pi^\rho_n(\xb)=\Pr_{\pi^\rho}\(\|\wb_n\|>\lambda n^{1/2}\)
&\le \frac{1}{\lambda^2 n}\E_{\pi^\rho_n} \|\wb_n\|^2
=\frac{2}{\lambda^2}\E_\mu |x|^2,
\end{align*}
by the Chebyshev inequality.
This shows that (recalling that $\pi^1=\mu\times \mu$)
\begin{equation}\label{Eq:outside}
\sum_{\|\xb\|>\lambda n^{1/2}}|\pi^\rho_n(\xb)-\mu_n\times \mu_n(\xb)|\le \frac{4}{\lambda^2}\E_\mu |x|^2.
\end{equation}
Since $\pi^\rho$ is aperiodic for $\rho \in (0, 1]$, 
then the local central limit theorem \cite[Theorem 2.3.9]{lawler_limic_2010} shows the following:
There exists a sequence $\delta_n>0$ such that $\delta_n \to 0$ as $n\to \infty$, for all $n \in \Z_{>0}$ and for all $\xb\in \Z^{2m}$,
\begin{equation}\label{Eq:LCLTZ}
|\pi^\rho_n(\xb)-f_{n\SS^\rho}(\xb)| \le \frac{\delta_n}{n^m} \quad \text{and} \quad |\mu_n\times \mu_n(\xb)-f_{n\SS^1}(\xb)|\le \frac{\delta_n}{n^m}.
\end{equation}
Note that there exists a constant $C_m>0$ such that for all real $\lambda>0$ and all $n \in \Z_{>0}$, the number of $\xb \in \Z^{2m}$ with $\|\xb\|\le \lambda n^{1/2}$ is at most $C_m \lambda^{2m} n^m$.
Therefore it holds that
\begin{align}\label{Eq:inside}
\sum_{\|\xb\|\le \lambda n^{1/2}}|\pi^\rho_n(\xb)-\mu_n\times \mu_n(\xb)| 
\le C_m\lambda^{2m} n^m \cdot\frac{2\delta_n}{n^m}+ \sum_{\|\xb\|\le \lambda n^{1/2}}|f_{n\SS^\rho}(\xb)-f_{n\SS^1}(\xb)|.
\end{align}
Hence by \eqref{Eq:limsup_rho}, \eqref{Eq:outside} and \eqref{Eq:inside}, for every $\lambda>0$,
\[
\limsup_{n \to \infty} \|\pi^\rho_n-\mu_n\times \mu_n\|_1 \le   3\lambda^2\|(\SS^\rho)^{-1}-(\SS^1)^{-1}\|+\left|1-\sqrt{\frac{\det \SS^\rho}{\det \SS^1}}\right|+\frac{4}{\lambda^2}\E_\mu|x|^2.
\]
Since $\SS^\rho \to \SS^1$ as $\rho \to 1$, 
we obtain
\[
\limsup_{\rho \to 1}\limsup_{n \to \infty} \|\pi^\rho_n-\mu_n\times \mu_n\|_1 \le \frac{4}{\lambda^2}\E_\mu|x|^2.
\]
This holds for all $\lambda>0$, and thus
we obtain in the total variation distance
\[
\lim_{\rho \to 1}\limsup_{n \to \infty}\|\pi^\rho_n-\mu_n\times \mu_n\|_\TV =0.
\]
This shows the upper bound in the claim.

Next let us show the lower bound in the claim.
The general case reduces to the case when $m=1$ since a projection $\Z^m \to \Z$ (whence $\Z^{2m} \to \Z^2$) to a coordinate only decreases the total variation distance.
Let $\SS=\SS^{\rho}$.
For $\xb_i \in \R^m$, $i=1, 2$,
\begin{align*}
|\abr{\xb_1, \SS^{-1}\xb_1}-\abr{\xb_2, \SS^{-1}\xb_2}| 
&\le 2\int_0^1 |\abr{\xb_2-\xb_1, \SS^{-1}(\xb_1+t(\xb_2-\xb_1))}|\,dt\\
&\le 2\|\SS^{-1}\|\|\xb_1-\xb_2\| \max\{\|\xb_1\|, \|\xb_2\|\}.
\end{align*}
Thus, if $\|\xb_1-\xb_2\|_\infty \le 1$ (where $\|\cdot\|_\infty$ denotes the supremum norm) and $\|\xb_i\| \ge \sqrt{2}$, $i=1, 2$, 
then $\|\xb_1-\xb_2\| \le \|\xb_2\|$ and
\[
|\abr{\xb_1, \SS^{-1}\xb_1}-\abr{\xb_2, \SS^{-1}\xb_2}| \le 2\|\SS^{-1}\|\|\xb_1-\xb_2\| (\|\xb_2\|+\|\xb_1-\xb_2\|) \le 4\sqrt{2}\|\SS^{-1}\|\|\xb_2\|.
\]
This shows that for $\xb_i \in \R^m$ such that $\|\xb_1-\xb_2\|_\infty \le 1$ and $\|\xb_i\|\ge \sqrt{2}$, $i=1, 2$,
\[
f_{n\SS}(\xb_1) \le f_{n\SS}(\xb_2)e^{\frac{c_\SS}{n}\|\xb_2\|}, \quad \text{where $c_\SS:=2\sqrt{2}\|\SS^{-1}\|$}.
\]
Therefore we obtain
\begin{align*}
\sum_{\|\xb\| >\lambda n^{\frac{1}{2}}}f_{n\SS}(\xb) \le \int_{\|\xb\|\ge \lambda n^{\frac{1}{2}}-\sqrt{2}}\frac{1}{\sqrt{(2\pi)^2 \det (n\SS)}}e^{-\frac{1}{2n}\abr{\xb, \SS^{-1}\xb}+\frac{c_\SS}{n}\|\xb\|}\,d\xb.
\end{align*}
Since $\SS^{-1}$ is positive definite, there exists a constant $\alpha>0$ such that $\abr{\xb, \SS^{-1}\xb} \ge \alpha \|\xb\|^2$ for all $\xb \in \R^m$, and
\[
-\frac{1}{2n}\abr{\xb, \SS^{-1}\xb}+\frac{c_\SS}{n}\|\xb\|
\le -\frac{\alpha}{2n}\|\xb\|^2+\frac{c_\SS}{n}\|\xb\|
=-\frac{\alpha}{2n}\(\|\xb\|-\frac{c_\SS}{\alpha}\)^2+\frac{c_\SS^2}{2\alpha n}.
\]
This shows that the last integral is at most, by the change of variables $\xb \mapsto \sqrt{n}\xb$,
\begin{align*}
&\int_{\|\xb\| \ge \lambda n^{\frac{1}{2}}-\sqrt{2}}\frac{1}{\sqrt{(2\pi)^2 \det (n\SS)}}e^{-\frac{\alpha}{2n}\(\|\xb\|-\frac{c_\SS}{\alpha}\)^2+\frac{c_\SS^2}{2\alpha n}}\,d\xb\\
&=\frac{e^{\frac{c_\SS^2}{2\alpha n}}}{2\pi \sqrt{\det \SS}}\int_{\|\xb\|\ge (\lambda n^{\frac{1}{2}}-\sqrt{2})/\sqrt{n}}e^{-\frac{\alpha}{2}\(\|\xb\|-\frac{c_\SS}{\alpha \sqrt{n}}\)^2}\,d\xb
=\frac{e^{\frac{c_\SS^2}{2\alpha n}}}{\sqrt{\det \SS}}\int_{(\lambda n^{\frac{1}{2}}-\sqrt{2})/\sqrt{n}}^\infty e^{-\frac{\alpha}{2}\(s-\frac{c_\SS}{\alpha\sqrt{n}}\)^2}s\,ds.
\end{align*}
In the above, the last equality has used the polar coordinate.
There exists a constant $C_\SS>0$ such that 
$s+c_\SS/(\alpha\sqrt{n}) \le s e^{\frac{\alpha}{4}s^2}$
for all $s>C_\SS$ and all $n\ge 1$.
The last term above equals by the change of variables, for $R:=(\lambda n^{\frac{1}{2}}-\sqrt{2})/\sqrt{n}-c_\SS/(\alpha \sqrt{n})$,
\[
\frac{e^{\frac{c_\SS^2}{2\alpha n}}}{\sqrt{\det \SS}}\int_{R}^\infty e^{-\frac{\alpha}{2}s^2}\(s+\frac{c_\SS}{\alpha \sqrt{n}}\)\,ds
\le \frac{e^{\frac{c_\SS^2}{2\alpha n}}}{\sqrt{\det \SS}}\int_{R}^\infty e^{-\frac{\alpha}{4}s^2}s\,ds
=\frac{e^{\frac{c_\SS^2}{2\alpha n}}}{\sqrt{\det \SS}}\frac{2}{\alpha}e^{-\frac{\alpha}{4}R^2}.
\]
Note that for all $n> (c_\SS/\alpha)^2+2$,
\[
R=(\lambda n^{\frac{1}{2}}-\sqrt{2})/\sqrt{n}-c_\SS/(\alpha \sqrt{n}) \ge\lambda -\sqrt{2}/\sqrt{n}-1 \ge \lambda-2.
\]
Hence there exist constants $c_\rho, C_\rho>0$ such that for all $n>C_\rho$ and all $\lambda \ge C_\rho$,
\[
\sum_{\|\xb\| >\lambda n^{\frac{1}{2}}}f_{n\SS}(\xb) \le C_\rho e^{-c_\rho \lambda^2}.
\]

Therefore together with \eqref{Eq:LCLTZ}, 
there exist (possibly different) constants $c_\rho, C_\rho>0$ depending only on $\rho$
such that the following holds: For all $n>C_\rho$ and all $\lambda \ge C_\rho$, 
\begin{align*}
&\|\pi^\rho_n-\mu_n\times \mu_n\|_1
\ge \sum_{\|\xb\| \le \lambda n^{1/2}}|\pi_n^\rho(\xb)-\mu_n\times \mu_n(\xb)|\\
&\ge \sum_{\|\xb\| \le \lambda n^{1/2}}|f_{n\SS^\rho}(\xb)-f_{n\SS^1}(\xb)|-2C_m\lambda^{2m} \delta_n
\ge \|f_{n\SS^\rho}-f_{n\SS^1}\|_1-2C_\rho e^{-c_\rho \lambda^2}-2C_m\lambda^{2m} \delta_n.
\end{align*}
Thus by \eqref{Eq:liminf_rho}, letting $n \to \infty$, we obtain for all small enough $\rho>0$,
\[
\liminf_{n\to \infty}\|\pi^\rho_n-\mu_n\times \mu_n\|_1 \ge  2-4\rho^{\frac{1}{4}}-2C_\rho e^{-c_\rho \lambda^2}.
\]
In total variation distance, letting $\lambda \to \infty$ and then $\rho \to 0$ yields
\[
\lim_{\rho \to 0}\liminf_{n \to \infty}\|\pi^\rho_n-\mu_n\times \mu_n\|_\TV =1,
\]
as required.
\qed

\subsection*{Acknowledgments}
The author would like to thank Professor J\'er\'emie Brieussel for helpful comments on an earlier version of this paper and discussions, and Professors Benoit Collins, Naotaka Kajino, Shu Kanazawa, Tomoki Kawahira, Tomoyuki Shirai and Kouji Yano for useful comments and discussions.
He also thanks the anonymous referees for constructive comments and suggestions.
The author is partially supported by 
JSPS Grant-in-Aid for Scientific Research JP20K03602 and JP24K06711.

\bibliographystyle{alpha}
\bibliography{nss}

\begin{thebibliography}{HSC93}

\bibitem[AB08]{AbramenkoBrown}
Peter Abramenko and Kenneth~S. Brown.
\newblock {\em Buildings}, volume 248 of {\em Graduate Texts in Mathematics}.
\newblock Springer, New York, 2008.
\newblock Theory and applications.

\bibitem[BB23]{BenjaminiBrieussel}
Itai Benjamini and J\'{e}r\'{e}mie Brieussel.
\newblock Noise sensitivity of random walks on groups.
\newblock {\em ALEA Lat. Am. J. Probab. Math. Stat.}, 20(2):1139--1164, 2023.

\bibitem[HSC93]{HebischSaloffCoste93}
W.~Hebisch and L.~Saloff-Coste.
\newblock Gaussian estimates for {M}arkov chains and random walks on groups.
\newblock {\em Ann. Probab.}, 21(2):673--709, 1993.

\bibitem[Kal18]{KalaiICM2018}
Gil Kalai.
\newblock Three puzzles on mathematics, computation, and games.
\newblock In {\em Proceedings of the {I}nternational {C}ongress of
  {M}athematicians---{R}io de {J}aneiro 2018. {V}ol. {I}. {P}lenary lectures},
  pages 551--606. World Sci. Publ., Hackensack, NJ, 2018.

\bibitem[KS83]{KramliSzasz}
Andr\'{a}s Kr\'{a}mli and Domokos Sz\'{a}sz.
\newblock Random walks with internal degrees of freedom. {I}. {L}ocal limit
  theorems.
\newblock {\em Z. Wahrsch. Verw. Gebiete}, 63(1):85--95, 1983.

\bibitem[KS01]{KotaniSunadaStandard}
Motoko Kotani and Toshikazu Sunada.
\newblock Standard realizations of crystal lattices via harmonic maps.
\newblock {\em Trans. Amer. Math. Soc.}, 353(1):1--20, 2001.

\bibitem[KSS98]{KotaniShiraiSunada}
Motoko Kotani, Tomoyuki Shirai, and Toshikazu Sunada.
\newblock Asymptotic behavior of the transition probability of a random walk on
  an infinite graph.
\newblock {\em J. Funct. Anal.}, 159(2):664--689, 1998.

\bibitem[LL10]{lawler_limic_2010}
Gregory~F. Lawler and Vlada Limic.
\newblock {\em Random Walk: A Modern Introduction}.
\newblock Cambridge Studies in Advanced Mathematics. Cambridge University
  Press, 2010.

\bibitem[PS94]{PollicottSharpRates}
Mark Pollicott and Richard Sharp.
\newblock Rates of recurrence for {${\bf Z}^q$} and {${\bf R}^q$} extensions of
  subshifts of finite type.
\newblock {\em J. London Math. Soc. (2)}, 49(2):401--416, 1994.

\bibitem[Sun13]{SunadaTopological}
Toshikazu Sunada.
\newblock {\em Topological Crystallography: With a View Towards Discrete
  Geometric Analysis}, volume~6 of {\em Surveys and Tutorials in the Applied
  Mathematical Sciences}.
\newblock Springer, Tokyo, 2013.

\bibitem[Tan24]{NonNS}
Ryokichi Tanaka.
\newblock Non-noise sensitivity for word hyperbolic groups.
\newblock {\em Ann. Fac. Sci. Toulouse Math. (6)}, 33(5):1487--1510, 2024.

\bibitem[Woe00]{WoessBook}
Wolfgang Woess.
\newblock {\em Random walks on infinite graphs and groups}, volume 138 of {\em
  Cambridge Tracts in Mathematics}.
\newblock Cambridge University Press, Cambridge, 2000.

\end{thebibliography}

\end{document}